\documentclass{article}

\usepackage{graphicx}

\usepackage[dvipsnames,usenames]{color}
\usepackage{caption}
\usepackage{amsmath, amssymb, amsthm, enumerate}
\usepackage{natbib}
\usepackage{bm}
\usepackage{changepage}

\usepackage{amsmath,amssymb,amsthm,enumerate,mathrsfs}
\newtheorem{theorem}{Theorem}[section]
\newtheorem{corollary}[theorem]{Corollary}
\newtheorem{lemma}[theorem]{Lemma}
\newtheorem{proposition}[theorem]{Proposition}
\theoremstyle{definition}
\newtheorem{definition}[theorem]{Definition}
\newtheorem{example}[theorem]{Example}

\newtheorem{assumption}{Assumption}[section]

\newcommand{\X}{{\cal X}}
\newcommand{\vt}{{\vartheta}}
\newcommand{\inmat}[1]{\mbox{\rm {#1}}}
\newcommand{\R}{{\rm I\!R}}
\newcommand{\bgeqn}{\begin{eqnarray}}
\newcommand{\edeqn}{\end{eqnarray}}
\newcommand{\bgeq}{\begin{eqnarray*}}
\newcommand{\edeq}{\end{eqnarray*}}
\newcommand{\dd}{\mathsf {d\kern -0.07em l}} 

\usepackage{algorithm} 
\usepackage{algorithmic}

\numberwithin{equation}{section}

\begin{document}
\makeatletter

\allowdisplaybreaks

\begin{center}
\large{\bf Bayesian Distributionally Robust Nash Equilibrium and Its Application}\footnote{ This paper is
dedicated to Professor R. T. Rockafellar on the occasion of his 90th birthday.
The work is supported by CUHK Start-up grant.
} 
\end{center}\vspace{5mm}
\begin{center}
\textsc{Jian Liu, Ziheng Su, and Huifu Xu}
\footnote{Department of Systems Engineering and Engineering Management.
The Chinese University of Hong Kong, Shatin, N.T., Hong Kong.
E-mail address: liujian@se.cuhk.edu.hk, zsu@se.cuhk.edu.hk, hfxu@se.cuhk.edu.hk.}
\end{center}

\vspace{2mm}

\footnotesize{
\noindent\begin{minipage}{14cm}
{\bf Abstract:}
Inspired by the recent work
of Shapiro et al.~\cite{shapiro2023bayesian},
we propose a Bayesian distributionally robust Nash equilibrium (BDRNE) model where each player lacks complete information on the true probability distribution of the underlying exogenous uncertainty represented by a random variable and subsequently 
determines the optimal decision by solving a Bayesian distributionally robust optimization (BDRO) problem under the Nash conjecture. Unlike most of the DRO models in the literature, the BDRO model assumes (a) the true unknown distribution of the random variable can be approximated by a randomized parametric family of distributions, (b) the average of the worst-case expected value of the objective function with respect to the posterior distribution of the parameter, instead of the worst-case expected value of the objective function is considered in each player's decision making, and (c) the posterior distribution of the parameter is updated as more and more sampling information of the random variable is gathered. Under some moderate conditions, we demonstrate the existence of a BDRNE and derive asymptotic convergence of the equilibrium as 
the sample size increases. Moreover, we propose to solve the BDRNE problem 
by Gauss-Seidel-type iterative method in the case when the ambiguity set of each player is constructed via Kullback-Leibler (KL) divergence. 
Finally, we apply the BDRNE model to a price competition problem
under multinomial logit demand. The preliminary numerical test results
show that the proposed model and computational scheme
perform well.
\end{minipage}
 \\[5mm]


\hbox to14cm{\hrulefill}\par


\section{Introduction}

Consider a 
%
game where $n$ players 
compete to provide goods or services in a non-collaborative manner. Let $u_j(x_j, x_{-j}): \R^{\sum_{j=1}^n n_j} \rightarrow \R$ 
be a continuous real-valued function 
representing 
player $j$'s payoff or its utility, with $x_j\in \mathcal{X}_j \subset \R^{n_j}$ being player $j$'s decision vector and $x_{-j}\in \mathcal{X}_{-j}$ being the decision variables of all players except player $j$. Assuming that each player makes a decision by maximizing 
his/her payoff/utility by holding the rival's strategies as fixed,
one may consider a situation where no player can be better off by unilaterally 
changing its position.
Such a situation is known as 
Nash 
equilibrium (NE) which 
is 
characterized by 
an $n$-tuple of the 
decision vectors of all players, written
$x^\ast:= (x_1^\ast, \cdots, x_n^\ast)\in \mathcal{X}:= (\mathcal{X}_1, \cdots, \mathcal{X}_n)$,
which satisfies 
$$
x_j^\ast \in \mathop{\arg \max}_{x_j\in \mathcal{X}_j}\ u_j(x_j, x_{-j}^\ast), 
\quad \inmat{for}\; j=1,\cdots,n,
$$   
where ${\cal X}_j$ denotes player $j$'s feasible set
of the decision vector $x_j$.
Since 
the foundational work 
by Nash
\cite{nash1950equilibrium, nash1950non},
Nash equilibrium  has been 
a 
widely used solution concept to describe the outcome of non-cooperative games 
in economics \citep{myerson1999nash}
and political science \citep{morrow1994game}.
Over the past few decades, various extensions of NE have been proposed driven by extensive practical applications
see e.g.,~\cite{facchinei2003finite, facchinei2010generalized, rockafellar2024generalized} and references therein for a comprehensive overview. 
In many real-world application problems 
such as market competition, transportation \cite{watling2006user,li2023modified}, and energy/signal network planning and engineering design,
decisions are often made under uncertainty 
due to limited observability of data, imprecise measurements, 
implementation variances,
and prediction inaccuracies. 
Consequently, one may consider 
a stochastic Nash equilibrium
(SNE) where each player determines the optimal response/decision by solving a stochastic optimization problem under the Nash conjecture: 
$$
x_j^\ast \in \mathop{\arg \max}_{x_j\in \mathcal{X}_j}\ \mathbb{E}_{Q}\big[u_j(x_j, x_{-j}^\ast, \xi)\big],  \quad \inmat{for}\; j=1,\cdots,n,
$$
where $\xi:(\Omega, {\cal F}, Q)\to \R^l$ 
 represents the underlying uncertainty
 and the mathematical expectation 
 is taken with respect to 
the probability distribution 
$Q$ of $\xi$. In this setup, each player 
aims to make an 
optimal decision before 
observing realization of the underlying uncertainty.
The SNE model 
is well studied,  see, e.g., \cite{jiang2008stochastic,xu2010sample,ravat2011characterization, xu2013stochastic, lei2022stochastic} and the references therein.

In some data-driven problems,
the true probability distribution 
$Q$ is often unknown and has to be 
approximated using partially available information such as sampling, subjective judgement based 
on empirical statistical data (e.g.~mean values and variances). 
Instead of using a single approximate 
probability distribution in classical stochastic programming models,
one may use the partially available information to construct a set of plausible probability distributions 
which has a higher likelihood to approximate and/or containing the true
and base the optimal decision on the worst-case probability distribution.
This kind of approach is known as distributionally robust optimization (DRO), which is aimed to mitigate the modeling risk arising from incomplete information of the true probability distribution, see \cite{rahimian2019distributionally} for an overview. Using the DRO paradigm, we may consider the following distributionally robust Nash equilibrium (DRNE) model,
which is an $n$-tuple 
$x^\ast:= (x_1^\ast, \cdots, x_n^\ast)\in \mathcal{X}:= (\mathcal{X}_1, \cdots, \mathcal{X}_n)$ satisfying 
\bgeqn 
\label{eq:DR-NE}
x_j^\ast \in \mathop{\arg \max}_{x_j\in \mathcal{X}_j} \min_{Q\in \mathcal{Q}_j} \mathbb{E}_Q \big[u_j(x_j, x_{-j}^\ast, \xi)\big],  \quad \inmat{for}\; j=1,\cdots,n,
\edeqn 
where $\mathcal{Q}_j$ denotes the ambiguity set of probability distributions constructed 
by player $j$. Since different players may have different sources to acquire 
information about the true probability distribution of $\xi$, the ambiguity sets 
are different unless there is a specific information structure 
which allows all players to share
the information. 
DRNE has been studied in several literature.
Loizou~\cite{loizou2016distributionally} proposes a DRNE model with each player’s objective being a conditional value-at-risk. One of the main focuses of Loizou’s work is to investigate cases where distributionally robust games are equivalent to Nash games without private information.
Sun et al.~\cite{sun2016convergence} present a convergence analysis of DRNE when players gain more and more information to construct a set of probability distributions as the ambiguity set.
Liu et al.~\cite{liu2018distributionally} develop several DRNE models, 
where some or all of the players in the games lack complete information on the true probability distribution of underlying uncertainty but they need to make a decision prior to the realization of such uncertainty.
More recently, Hori and Yamashita \cite{hori2022two}
consider a two-stage distributionally robust noncooperative game 
with application to
Cournot–Nash competition
and demonstrate
existence of Nash equilibria. 
In the case that $\mathcal{Q}_j$ contains all of the probability 
measures/distributions in $\R^l$, (\ref{eq:DR-NE}) reduces to distribution-free robust Nash equilibrium \cite{aghassi2006robust}.

In this paper, we revisit the DRNE model by adopting the
Bayesian distributionally robust optimization (BDRO)
approach recently 
proposed by Shapiro et al.~\cite{shapiro2023bayesian}.
Specifically,
we consider a situation where the true unknown probability distribution of $\xi$ can be represented or approximated by a family of probability distributions 
$\{Q_\theta\}$ parameterized by some $\theta \in \Theta\subset \R^{t}$,
each player is aware of this and uses sample information 
to acquire a posterior 
distribution
of $\theta$ by virtue of
Bayesian formula
\begin{equation}
\label{eq:posterior distribution}
\rho(\theta | \bm{\xi}^{(N)}) = \frac{f(\bm{\xi}^{(N)} | \theta)\rho(\theta)}{\int_\Theta f(\bm{\xi}^{(N)} | \theta) \rho(\theta) d\theta},   
\end{equation}
where $\rho(\theta)$ and 
$\rho(\theta | \bm{\xi}^{(N)})$  denote
the prior and posterior density functions 
of $\theta$ respectively, $\bm{\xi}^{(N)}:=\{\xi_1, \cdots, \xi_N\}$ are independent and identically distributed (i.i.d.) samples, and
$f(\bm{\xi}^{(N)} | \theta) = \prod_{i=1}^N f(\xi_i | \theta)$ denotes the conditional density of the sample. 
Let $\mu_\theta$ and $\mu_{\theta}^N$ denote the 
respective prior and posterior distributions.
Subsequently, each player, say player $j$, $j=1,\cdots,n$,
determines its optimal decision by solving a BDRO problem
\begin{equation}\label{eq:BDRO}
\max_{x_j\in \mathcal{X}_j} \mathbb{E}_{\theta_j^{N_j}} 
\left[ \inf_{Q\in \mathcal{Q}_{\epsilon_j}^{\theta_j}} \mathbb{E}_{Q}\big[ u_j(x_j, x_{-j}, \xi) 
\big] \right],  
\end{equation}
where the expectation $\mathbb{E}_{\theta_j^{N_j}}$ 
is taken with respect to the
posterior distribution $\mu_{\theta_j}^{N_j}$ with density
function
$\rho(\theta_j| \bm{\xi}^{(N_j)})$
whereas the expectation $\mathbb{E}_{Q}$ is taken with respect to 
the distribution $Q$ from the ambiguity set constructed via Kullback-Leibler divergence \cite{kullback1951information} (KL-divergence) 
\begin{equation}\label{eq:ambiguitysetQ_j}
\mathcal{Q}_{\epsilon_j}^{\theta_j} := 
\big\{ Q: \dd_{KL}\big(Q || Q_ {\theta_j}\big) \leq \epsilon_j \big\}, \quad \inmat{for}\; \theta_j \in \Theta_j.
\end{equation}
The KL-divergence
$$
\dd_{KL} \big(Q||Q_{\theta_j}\big) := \int_\Xi q(\xi) \ln\big(q(\xi) / f_{\theta_j}(\xi) \big) d\xi
$$ 
measures the deviation 
of 
distribution $Q$ 
from the nominal distribution $Q_{\theta_j}$ via their density functions 
$q(\xi)$ and $f_{\theta_j}(\xi):= f(\xi| \theta_j)$.
The KL-divergence is a special case of the $\phi$-divergence and has been studied
in DRO, see the earlier work by Hu and Hong \cite{hu2013kullback}.
If we interpret the quantity as a 
``distance'', 
then the ambiguity set has a kind of ``ball'' structure centered at the nominal distribution $Q_{\theta_j}$
with “radius” $\epsilon_j$. The latter signifies 
the decision maker's confidence in the ambiguity
set containing the true unknown probability distribution.
In the case that $\xi$ is discretely distributed,
$q(\xi)$ and $f_{\theta_j}(\xi)$ may be viewed as probability mass functions.
Since the ambiguity set defined by the KL-divergence is compact under the topology of weak convergence \cite[Theorem 20]{van2014renyi},
the infimum in (\ref{eq:BDRO}) is attainable. 
The attainability of the outer maximum requires some more sophisticated arguments, we will come back to this in Section \ref{se:BDRNE}.
Assuming that each player determines the
optimal decisions by solving (\ref{eq:BDRO}),
we may consider 
an equilibrium
which we call Bayesian distributionally robust Nash equilibrium (BDRNE).

\begin{definition}[BDRNE]
\label{def:BDRNE} 
The BDRNE is an $n$-tuple $\hat{x}^N_\epsilon := (\hat{x}_{\epsilon_1}^{N_1}, \cdots, \hat{x}_{\epsilon_n}^{N_n}) \in \mathcal{X}:= (\mathcal{X}_1, \cdots, \mathcal{X}_n)$ such that
\begin{equation}
\label{eq:BDRNE}
\hat{x}_{\epsilon_j}^{N_j} \in \mathop{\arg \max}_{x_j\in \mathcal{X}_j}\ \mathbb{E}_{\theta_j^{N_j}} 
\left[ \inf_{Q\in \mathcal{Q}_{\epsilon_j}^{\theta_j}} \mathbb{E}_{Q}\big[ u_j(x_j, \hat{x}_{\epsilon_{-j}}^{N_{-j}}, \xi) \big] \right],   \quad \inmat{for}\; j=1,\cdots,n,
\end{equation}
where 
the superscript $N_j$ and subscript $\epsilon_j$ in 
$\hat{x}_{\epsilon_j}^{N_j}$ 
signify that the equilibrium depends on both the ``radius'' 
$\epsilon_j$ of 
the ambiguity set and 
the sample size $N_j$.
\end{definition}

In the BDRNE model, each player determines its best response strategy under the Nash conjecture, assuming that its rivals' best response strategy $\hat{x}_{\epsilon_{-j}}^{N_{-j}}$ is fixed, where $N_{-j}$ and $\epsilon_{-j}$ represent the vector of all players except player $j$, with $N:= (N_1, \cdots, N_n)$ and $\epsilon:= (\epsilon_1, \cdots, \epsilon_n)$.  
It means that player $j$ does not have to know the rivals' optimization problems. 
The analysis of 
existence and uniqueness of an equilibrium 
in the forthcoming discussions will be based on this conjecture.
We can also interpret the model
from a modeler's perspective.
    Assuming that each player is risk adverse against the risk arising from incomplete probability distribution of exogenous uncertainty 
    $\xi$ and each uses the BDRO model for optimal decision-making under the Nash conjecture, the modeller would like to understand the 
    outcome of an equilibrium
    from such a game.
    The model may provide
    a mathematical 
    approach for analyzing 
    market equilibrium such as price competition to be discussion in Section 3. It does not mean that in practice 
    every player actually 
    uses KL-divergence to construct an ambiguity set 
    as presented in the formulation
    and then works out an optimal strategy.
    Rather it means that if each player does as 
    formulated above, then there will be an equilibrium under some moderate conditions and such as equilibrium 
   which  may be used to explain/understand market competition or prediction of outcome of market competition in future.

The ``radius''
$\epsilon_j$ depends on player $j$'s confidence on the reliability of the nominal distribution $Q_{\theta_j}$ with density function $f_{\theta_j}(\xi)$. Unlike the DRNE models in the literature where the ambiguity set is constructed via statistical quantities such as mean and covariance,
here we consider a different information structure which each player may elicit the true probability distribution of $\xi$:
it allows each player to draw samples of $\xi$
to 
update the posterior 
distribution of the random parameter. 
As the simple sizes increases,
the posterior 
distribution converges to some distribution. 
In that case we may set 
$\epsilon_j = 0$, and 
$\mathcal{Q}_{\epsilon_j}^{\theta_j}$ reduces to a singleton containing only $Q_{\theta_j}$, 
consequently player $j$'s
BDRO problem \eqref{eq:BDRO} reduces to the Bayesian average 
(BAVG) problem
\begin{equation}\label{eq:BAVG}
\max_{x_j\in \mathcal{X}_j} \mathbb{E}_{\theta_j^{N_j}} 
\left[ \mathbb{E}_{\xi | \theta_j} \big[ u_j(x_j, x_{-j}, \xi) \big] \right],
\end{equation}
where the expectation $\mathbb{E}_{\xi | \theta_j}$ is taken with respect to the nominal distribution $Q_{\theta_j}$ with density function $f_{\theta_j}(\xi)$ 
and
the expectation $\mathbb{E}_{\theta_j^{N_j}}$
is taken with respect to the
posterior distribution $\mu_{\theta_j}^{N_j}$ with density
function
$\rho(\theta_j| \bm{\xi}^{(N_j)})$.
In that case,
each player uses the nominal distribution in 
optimal decision-making.
We call the outcome of this kind of competition 
the Bayesian average Nash equilibrium (BANE) which is formally defined as follows.
\begin{definition}[BANE]\label{def:Bayesian-Average-NE} 
The BANE is an $n$-tuple $\bar{x}^N:= (\bar{x}_1^{N_1}, \cdots, \bar{x}_n^{N_n})\in \mathcal{X}:= (\mathcal{X}_1, \cdots, \mathcal{X}_n)$ such that
\begin{equation}\label{eq:Bayesian-Average}
\bar{x}_j^{N_j} \in \mathop{\arg \max}_{x_j\in \mathcal{X}_j}\ \mathbb{E}_{\theta_j^{N_j}} 
\left[ \mathbb{E}_{\xi| \theta_j}\big[ u_j(x_j, \bar{x}^{N_{-j}}_{-j}, \xi) \big] \right],   \quad \inmat{for}\; j=1,\cdots,n.  
\end{equation}
\end{definition}

In the BANE model, the ``average'' means that 
$\mathbb{E}_{\theta_j^{N_j}} \big[\mathbb{E}_{\xi| \theta_j}[\cdot]\big]
$ is the expected value w.r.t. the joint distribution of $\xi$ and $\theta_j$. Note that an unbiased estimate of $\mathbb{E}_{\theta_j^{N_j}} \big[\mathbb{E}_{\xi| \theta_j}[\cdot]\big]$ can be obtained by generating a random sample of $\theta_j$ from the 
posterior distribution $\mu_{\theta_j}^{N_j}$ with density $\rho(\theta_j| \bm{\xi}^{(N_j)})$ and then generating a random realization of $\xi \sim Q_{\theta_j}$ conditional on the generated $\theta_j$. This enables us to apply either the sample average approximation (SAA) or the stochastic approximation optimization method to solve problem \eqref{eq:Bayesian-Average}, provided that there is an efficient way to generate such random samples.
At this point, it might be helpful to 
distinguish BANE from the well-known Bayesian Nash equilibrium (BNE) \cite{harsanyi1995games}. The latter considers 
each player's utility function being characterized by a type parameter 
and information about the type parameter is private information. By randomizing rivals' type parameters,
and assuming their optimal decisions are dependent on the parameters, each player sets  his/her optimal decision by maximizing the expected value of his/her utility  where the expectation is taken with respect to the joint probability distribution of all type 
the rivals' parameters. 
Thus, 
BANE is concerned with exogenous uncertainty whereas
BNE is
about 
endogenous (DMs') uncertainty
related to players'
utility/risk preferences. 
From mathematical structure perspective, the former is 
one-stage stochastic program (decision is independent of the underlying uncertainty) whereas the latter is two-stage (because the decision is uncertainty dependent), see \cite{meirowitz2003existence, guo2021existence, liusunxu2024Bayesian, tao2025generalized, su2025continuous} for the recent development of the latter.



The main contributions of this paper can be summarized as follows.
\begin{itemize}
 
\item  \textbf{Modelling}.
We propose a one-stage 
BDRNE model
where each player solves a BDRO problem under the Nash conjecture.
Unlike the DRNE model proposed by Liu et al.~\cite{liu2018distributionally} where the ambiguity set of probability distributions in
each player's decision-making problem 
is constructed under moment-type information structure, 
individual player's ambiguity set in the BDRNE model is constructed 
with a ``ball'' structure centered 
at a nominal distribution using KL-divergence.
Moreover,  
the true unknown probability distribution is
approximated by a family of parametric distributions, and
the posterior distribution of the parameter is
updated via the Bayesian 
approach as 
more and more 
samples of the underlying uncertainty 
are
drawn 
by the player.
Furthermore,
since the ambiguity set is affected by each realization of the parameter via its nominal 
distribution,
the average of the worst-case expected value 
of the objective function with respect to the posterior distribution of the parameter, instead of the worst-case expected value 
of the objective function is considered in each player's decision-making.

\item  \textbf{Theory and numerical methods}. 
Under some moderate conditions, we derive 
the existence of a BDRNE
and 
demonstrate that such an equilibrium is stable against data perturbation. Moreover,
since the equilibrium is ad hoc 
with available sample information, we investigate 
asymptotic convergence of 
BDRNE as sample size increases.
To solve the equilibrium problem,
we follow Shapiro et al.~\cite{shapiro2023bayesian} to
reformulate each player's maximin problem as a 
single 
one-stage stochastic maximization problem and discretize it via  
SAA.
We then solve 
the latter 
by the Gauss-Seidel-type iteration method~\cite{facchinei2010generalized}.
This differs from Liu et al.~\cite{liu2018distributionally} which solves each player's minimax problem in an alternating manner.

\item \textbf{Application}.
We apply the proposed BDRNE model 
in 
a price competition problem 
under the multinomial logit (MNL) demand (BDRNE-MNL). 
This includes the derivation of conditions for ensuring the existence of an equilibrium,
stability of 
the equilibrium against 
variation of individual firm's 
marginal cost and finally numerical solution of an approximate equilibrium in price.

\end{itemize}

The rest of the paper is organized as follows. 
In Section \ref{se:BDRNE}, we 
analyze the BDRNE model
by establishing 
the existence, stability, and convergence of the BDRNE. 
Moreover, we propose 
an iterative scheme 
to compute an approximate BDRNE.
In Section \ref{se:MNL}, we 
apply the BDRNE model to
a price competition problem 
under the multinomial logit demand.
In Section \ref{se:Numerical Results}, we report the result of the numerical tests on the proposed model and computational scheme. 
Finally, in Section \ref{se:Concluding remarks}, we conclude the paper with remarks highlighting 
potential further work for future research. 

\section{BDRNE}\label{se:BDRNE}
In this section,
we discuss in detail the BDRNE model including 
well-definedness of the model, 
existence of an equilibrium, stability of the equilibrium, and asymptotic convergence of the equilibrium as the sample size increases.

\subsection{Existence of BDRNE}

We begin by discussing well-definedness and existence of an equilibrium.
To this end, 
we need to 
make the following assumption.

\begin{assumption}\label{ass:Assumption_1}
Let $u_j(x_j, x_{-j}, \xi)$, $j = 1, \cdots, n$,
be defined as in \eqref{eq:BDRO}. The following conditions are satisfied.
\begin{itemize}
    \item [(a)] $u_j(x_j, x_{-j}, \xi): \R^{\sum_{j=1}^n n_j} \times \R^l \rightarrow \R$ is continuous 
    and for each fixed $(x_{-j}, \xi)$, $u_j(\cdot, x_{-j}, \xi)$ is concave over $\mathcal{X}_j$.
    \item [(b)] The feasible sets $\mathcal{X}_j$, for $j = 1, \cdots, n,$ are compact. 
     \item [(c)] There exists a 
     non-negative integrable function 
     $\eta_j(\xi)$ such that $\mathbb{E}_{Q}[\eta_j(\xi)]<\infty$, and
     \bgeq
     \big|u_j(x_j, x_{-j}, \xi)\big|\leq \eta_j(\xi), \quad \forall x_j\in {\cal X}_j, x_{-j}\in \mathcal{X}_{-j}, \xi\in \Xi.
     \edeq 
     
     \item[(d)] 
     Let $\mathscr{P}(\R^l)$ denote the set of probability measures defined over $\R^l$ 
     and 
      \bgeqn 
      \label{eq:M_eta_j-PJOR}
    \mathcal{M}_{\R^l}^{\eta_j} :=\big\{ Q \in \mathscr{P}(\R^l): \mathbb{E}_{Q}\big[\max(\eta_j(\xi),\|\xi\|_\infty)\big] < \infty \big\}.
    \edeqn
    For $j=1,\cdots, n$, $Q_ {\theta_j}\in \mathcal{M}_{\R^l}^{\eta_j}$.
    Let
    \bgeq\label{eq:ambiguitysetQ_j-1}
    \mathcal{Q}_{\epsilon_j}^{\theta_j} := 
    \big\{ Q\in \mathcal{M}_{\R^l}^\eta: \dd_{KL}\big(Q || Q_ {\theta_j}\big) \leq \epsilon_j \big\}, \quad \inmat{for}\; \theta_j \in \Theta_j.
    \edeq
    \item[(e)] The density function $f_{\theta_j}$ is continuous in $\theta_j$, for $j=1,\cdots,n$.


\end{itemize}
\end{assumption}

Assumption~\ref{ass:Assumption_1} (a) and (b) is standard, see e.g.~\cite{liu2018distributionally}.
The compactness of ${\cal X}_j$ may be replaced 
by supcompactness condition (upper level set being compact). We impose the compactness 
condition to ease the discussion as our focus is not on the issue. Assumption~\ref{ass:Assumption_1} (c) is a kind of growth condition where $\eta_j(\xi)$ is not necessarily bounded. Assumption~\ref{ass:Assumption_1} (d)
restricts the probability distributions 
in the ambiguity sets to a subset of $\mathscr{P}(\R^l)$ with finite first order moment (mean value) of $\eta_j(\xi)$. It has implications on the data structure in the subsequent analysis. In the case that $\eta_j(\xi)$ is bounded,  $\mathcal{M}_{\R^l}^\eta =\mathscr{P}(\R^l)$. 
Assumption \ref{ass:Assumption_1} (e) is also satisfied by many parametric distributions, see e.g.~\cite{shapiro2023bayesian}.
Under Assumption~\ref{ass:Assumption_1}, each player's decision-making problem \eqref{eq:BDRNE} is well defined. The next proposition addresses this.

\begin{proposition}[Well-definedness of 
(\ref{eq:BDRNE})]
\label{pro:pro1} 
For $j=1,\cdots,n$,
let
\bgeq
\label{eq:BDRNE-v_j}
v_j(x_j,x_{-j},\theta_j):=\inf_{Q\in \mathcal{Q}_{\epsilon_j}^{\theta_j}} \mathbb{E}_{Q}\big[ u_j(x_j, 
x_{-j}, \xi) \big], \quad \inmat{and} \quad
\vt_j(x_j,x_{-j}):= \mathbb{E}_{\theta_j^{N_j}} 
 \big[ v_j(x_j,x_{-j},\theta_j)
\big].
\edeq
Under Assumption ~\ref{ass:Assumption_1},
the following assertions hold.
 
\begin{itemize}
    \item[(i)] The ambiguity set $\mathcal{Q}_{\epsilon_j}^{\theta_j}$ is
    continuous in $\theta_j$ under the topology of weak convergence.
    
    \item[(ii)]$v_j(x_j,x_{-j},\theta_j)$ is continuous in $(x_j,\theta_j)$. Moreover, if $\eta_j(\xi)$ is bounded, then
    \bgeqn
    \big|v_j(x_j,x_{-j},\theta_j)\big|\leq 
    \mathbb{E}_{Q}\big[e^{\eta_j(\xi)}\big]
    +\epsilon_j, \quad \inmat{for}\; Q\in \mathcal{Q}_{\epsilon_j}^{\theta_j}.
    \label{eq:v_j-bound}
    \edeqn

    \item[(iii)] If, in addition, 
    $\mathbb{E}_{\theta_j^{N_j}} \big[\sup_{x_j\in {\cal X}_j} v_j(x_j,x_{-j},\theta_j)\big]<\infty$, then 
    the maximum of problem \eqref{eq:BDRNE} is attainable and 
    the optimal value function $\vt_j(x_j)$ is continuous over ${\cal X}_j$.

\item[(iv)] Assume further that
there exist
positive constants $\varpi$ and  $\varrho$ such that
\bgeqn 
\label{eq:u_j-Lip-xi-Uinform-x}
\big|u_j(x_j,x_{-j},\xi')-
u_j(x_j,x_{-j},\xi)\big|
\leq \varpi \|\xi'-\xi\|, \quad \forall \xi', \xi\in\Xi,
\label{eq:U_j-Lip}
\edeqn 
and
\bgeqn 
\label{eq:ln-f_theta-Lip}
\big|\ln f_{\theta_j'}(\xi)-\ln f_{\theta_j}(\xi)\big|\leq \varrho\|\theta_j'-\theta_j\|, \quad  \forall \theta_j', \theta_j\in\Theta_j.
\label{eq:U_j-Lip-theta_j},
\edeqn 
and that $\Theta_j$ is a compact set.
Then
there exists a positive constant $C_1$ such that
\bgeqn
\sup_{x_j\in {\cal X}_j} \big|v_j(x_j,x_{-j},\theta_j')-
v_j(x_j,x_{-j},\theta_j)\big|\leq C_1\|\theta_j'-\theta_j\|, \quad \forall \theta_j',\theta_j \in\Theta_j.
\label{eq:uniform-v_j-prop}
\edeqn
 \end{itemize}
\end{proposition}

Before providing a proof, it might be helpful to make 
a few comments on some of the conditions and conclusions of the proposition.
Part (i) means that the ambiguity set $\mathcal{Q}_{\epsilon_j}^{\theta_j}$
changes continuously as $\theta_j$ varies. It paves
the way for 
deriving 
continuity of $v_j$ in $\theta_j$ in Part (ii).
The latter ensures that $v_j$ is a continuous random function of $\theta_j$ and hence its mathematical expectation 
w.r.t.~the distribution of $\theta_j$ is well-defined.
Condition (\ref{eq:u_j-Lip-xi-Uinform-x})
requires $u_j$ to be Lipschitz continuous over $\Xi$ uniformly for all $x\in {\cal X}$ and condition (\ref{eq:ln-f_theta-Lip}) requires
$f_{\theta_j}(\xi)$ to be Lipschitz continuous in $\theta_j$.
They may be satisfied for some specific utility functions and parametric 
family of distributions.
The two conditions are specifically imposed to derive (\ref{eq:uniform-v_j-prop}). In the case when $\Theta_j^*$ is a singleton, we no longer require these conditions to establish 
(\ref{eq:uniform-v_j-prop}). We will come back to this shortly.


\noindent
\textbf{Proof.} \underline{Part (i)}. 
Consider
the parametric system of inequality in \eqref{eq:ambiguitysetQ_j}.
Let
\bgeqn 
\label{eq:densty-KL}
\psi(q,\theta_j) :=\int_\Xi q(\xi) \ln\big(q(\xi) / f_{\theta_j}(\xi) \big) d\xi -\epsilon_j.
\edeqn 
Then $\psi(q,\theta_j)$ is convex in $q$
given that $x\ln x$ is a convex function and 
$
\psi(f_{\theta_j},\theta_j) = -\epsilon_j<0,
$
which means that 
the inequality system \eqref{eq:densty-KL}
satisfies the Slater constraint qualification 
uniformly for all $\theta_j$.
Let ${\cal H}_{\theta_j}$
denote the set of density functions
corresponding to the probability distributions 
in set $\mathcal{M}_{\R^l}^{\eta_j}$ which are
absolutely continuous with respect to $Q_{\theta_j}$.
Define:
\bgeq
{\cal F}_{\theta_j} :=\left\{q \in {\cal H}_{\theta_j}: 
\psi(q,\theta_j) \leq 0
\right\}.
\edeq
We want to show ${\cal F}_{\theta_j} $ is continuous in $\theta_j$. Based on the discussions above, we know that  
${\cal F}_{\theta_j}$ is a convex set and has nonempty interior.
For $q \in {\cal H}_{\theta_j}$, we define
$$
\hat{q} := \left(1-\frac{\upsilon}{\upsilon+\epsilon_j}\right)q+
\frac{\upsilon}{\upsilon+\epsilon_j}f_{\theta_j},
$$
where $\upsilon := |(\psi(q,\theta_j))_+|$, and $(a)_+ := \max(0,a)$, for $a\in \R$.
Then, we have
\begin{align*}
\psi(\hat{q},\theta_j) & \leq  \left(1-\frac{\upsilon}{\upsilon+\epsilon_j}\right) \psi(q,\theta_j) + 
 \frac{\upsilon}{\upsilon+\epsilon_j} \psi(f_{\theta_j},\theta_j)\leq  \left(1-\frac{\upsilon}{\upsilon+\epsilon_j}\right) 
\upsilon+ \frac{\upsilon}{\upsilon+\epsilon_j} 
(-\epsilon_j) = 0,
\end{align*}
which shows that $\hat{q}\in {\cal F}_{\theta_j}$, and hence
\bgeq
\dd(q,{\cal F}_{\theta_j}) \leq 
\dd(q,\hat{q}),
\edeq
where $\dd$ is the distance to be introduced in the space of density functions ${\cal H}_{\theta_j}$.
Next, we define the distance for the density functions.
For any $q_1,q_2\in {\cal H}_{\theta_j}$, we define
\bgeq
\label{eq:dist-density}
\dd(q_1,q_2) := \sup_{g\in {\cal G}}
\left|\int_\Xi g(\xi)q_1(\xi)d\xi -
\int_\Xi g(\xi)q_2(\xi)d\xi\right|, 
\edeq
where ${\cal G}$ is a class of functions
such that
$
|g(\xi_1)-g(\xi_2)|\leq \|\xi_1 - \xi_2\|,
$ for all $\xi_1,\xi_2\in\Xi$.
Since $\int_\Xi g(\xi)q_1(\xi)d\xi $ and 
$
\int_\Xi g(\xi)q_2(\xi)d\xi$
can be written as $\mathbb{E}_{Q_1}[g]$ and $\mathbb{E}_{Q_2}[g]$ respectively, where
$Q_1$ and $Q_2$ are probability distributions corresponding to $q_1$ and $q_2$, then
\bgeqn 
\label{eq:dist-density-1}
\dd(q_1,q_2) = \dd_K(Q_1,Q_2),
\edeqn 
where $\dd_K$ denotes the Kantorovich distance
of two probability distributions, i.e.,
\bgeqn 
\dd_K(Q_1,Q_2)=
\sup_{g\in \mathscr{G}}
\left| 
\mathbb{E}_{Q_1}\big[g\big] - \mathbb{E}_{Q_2}
\big[g\big]\right|,
\edeqn
where $\mathscr{G}$ denote the set of Lipschitz continuous functions over $\Xi$ with modulus $1$.
Note that under condition (\ref{eq:M_eta_j-PJOR}), $\int_\Xi g(\xi)q_i(\xi)d\xi$, for $i=1,2$, is well-defined.
It follows by \cite[Proposition 2.1]{zhang2024statistical}
that $\dd_K$ metrizes the topology of weak convergence in the space of probability measures
$\mathscr{P}(\R^l)$. Moreover,
by \cite[Theorem 20]{van2014renyi},
we know that $\mathcal{Q}_{\epsilon_j}^{\theta_j}$ is compact
 under the topology of weak convergence. Thus 
there exists a positive constant $C_2$ such that
\bgeqn
\label{eq:q-Q-theta_j-PJO}
\sup_{q_1,q_2 \in {\cal F}_{\theta_j}}\dd(q_1,q_2)
=
\sup_{Q_1,Q_2\in \mathcal{Q}_{\epsilon_j}^{\theta_j}}
\dd_K(Q_1,Q_2)\leq C_2.
\edeqn
Observe that
$
\epsilon_j \dd(q,\hat{q}) = \upsilon (\dd(q, f_{\theta_j}) - \dd(q,\hat{q})).
$
Thus
\bgeq
\dd(q,\hat{q}) =
\frac{\upsilon}{\epsilon_j}
(\dd(q, f_{\theta_j})  - \dd(q,\hat{q}))\leq 
\frac{\upsilon}{\epsilon_j} \dd(f_{\theta_j}, \hat{q})\leq
\frac{C_2}{\epsilon_j} |(\psi(q,\theta_j))_+|,
\edeq
and hence
\bgeq
\dd(q,{\cal F}_{\theta_j})\leq \dd(q,\hat{q})\leq 
\frac{C_2}{\epsilon_j} |(\psi(q,\theta_j))_+|.
\edeq
To show that 
${\cal F}_{\theta_j} $ is continuous in $\theta_j$, let us consider any fixed  $\theta_j\in \Theta_j$  and $\theta_j'\in \Theta_j$. Let
$q\in {\cal F}_{\theta_j'}$ be any density function in the set. Then 
$(\psi(q,\theta_j'))_+=0$, and hence
\bgeqn \dd(q,{\cal F}_{\theta_j}) \leq 
\frac{C_2}{\epsilon_j} \big| (\psi(q,\theta_j))_+ - 
(\psi(q,\theta_j'))_+ \big| \leq \frac{C_2}{\epsilon_j} 
\int_\Xi q(\xi) \big|\ln\big(f_{\theta_j'}(\xi)/ f_{\theta_j}(\xi) \big)\big| d\xi. 
\label{eq:q,F_theta_j}
\edeqn
The rhs of the inequality above goes to zero as $\theta_j'\to\theta_j$, which implies that
\bgeqn 
\label{eq:prop-2.1-parti-PJO-a}
\mathbb{D}({\cal F}_{\theta_j'},{\cal F}_{\theta_j}) :=\sup_{q\in {\cal F}_{\theta_j'}}
\dd(q, {\cal F}_{\theta_j})\to 0.
\edeqn 
This shows that
${\cal F}_{\theta_j}$ is upper semi-continuous at $\theta_j$ and hence  
$\mathcal{Q}_{\epsilon_j}^{\theta_j}$  is upper semi-continuous at $\theta_j$. 
By adjusting (increasing) $C_2$ if necessary, we assert that 
\eqref{eq:q-Q-theta_j-PJO} holds for all $\theta_j'$ near $\theta_j$. 
By swapping the positions between $\theta_j$ and $\theta_j'$ in the analysis above, we can establish 
\bgeqn 
\label{eq:prop-2.1-parti-PJO-b}
\mathbb{D}({\cal F}_{\theta_j},{\cal F}_{\theta_j'}) 
\to 0.
\edeqn 
Combining \eqref{eq:prop-2.1-parti-PJO-a}-\eqref{eq:prop-2.1-parti-PJO-b}, we assert that
${\cal F}_{\theta_j}$ is continuous in $\theta_j$
under the Hausdorff distance (of two compact sets 
in the space of 
${\cal H}_{\theta_j}$)
equipped with $\dd(\cdot)$.
By (\ref{eq:dist-density-1}), the continuity implies the continuity of $\mathcal{Q}_{\epsilon_j}^{\theta_j}$ in $\theta_j$.

\underline{Part (ii)}. 
Under 
Assumption~\ref{ass:Assumption_1}, 
$\mathbb{E}_{Q}\big[ u_j(x_j, \hat{x}_{\epsilon_{-j}}^{N_{-j}}, \xi) \big]$ is finite-valued and 
is continuous in $x_j$. Together with Part (i), 
we can use Berge's maximum theorem \cite{berge1877topological} to claim that 
the optimal value function $v_j(x_j,x_{-j},\theta_j)$ is continuous in both $x_j$ and $\theta_j$.
To show (\ref{eq:v_j-bound}), we note that by 
Donsker-Varadhan (DV) variational formula (see e.g.~\cite{birrell2022optimizing}),
\bgeq
\dd_{KL}(Q|| Q_{\theta_j})  = \sup_{h\in {\cal H}_b(\xi)}\big\{
\mathbb{E}_{Q}[h] - \log \mathbb{E}_{Q_{\theta_j}}[e^h]
  \big\},
\edeq
where ${\cal H}_b$ denotes the bounded functions
defined over $\R$. Under Assumption~\ref{ass:Assumption_1} (c), we have
$u_j(x_j, x_{-j}, \xi)\in {\cal H}_b(\xi)$ when $\eta_j(\xi)$ is bounded. Consequently, for $Q\in \mathcal{Q}_{\epsilon_j}^{\theta_j}$,
\bgeq 
\mathbb{E}_{Q}\big[u_j(x_j, x_{-j}, \xi)\big] \leq \log \mathbb{E}_{Q_{\theta_j
}}\big[e^{u_j(x_j, x_{-j}, \xi)}\big] +\epsilon_j\leq 
\log \mathbb{E}_{Q_{\theta_j}}\big[e^{\eta_j(\xi)}\big]
+\epsilon_j.
\edeq

\underline{Part (iii)}. By definition,
$\vt_j(x_j,x_{-j}) := \mathbb{E}_{\theta_j^{N_j}} [v_j(x_j,x_{-j},\theta_j)].$
The conclusion follows by Part (ii) and the Lebesgue dominated convergence theorem.

\underline{Part (iv)}.
Under conidition (\ref{eq:U_j-Lip}), we have
\bgeq
\left|\mathbb{E}_{Q'}\big[u_j(x_j, x_{-j}, \xi)\big]
-
\mathbb{E}_{Q}\big[u_j(x_j, x_{-j}, \xi)\big] \right|
\leq 
\varpi\dd_K(Q',Q).
\edeq
On the other hand, we can use conditions (\ref{eq:U_j-Lip-theta_j}) and (\ref{eq:q,F_theta_j}) 
to derive
\bgeq \dd(q,{\cal F}_{\theta_j}) 
\leq \frac{C_2}{\epsilon_j} 
\int_\Xi q(\xi) \big|\ln\big(f_{\theta_j'}(\xi)/ f_{\theta_j}(\xi) \big)\big| d\xi
\leq \frac{C_2}{\epsilon_j} \cdot \varrho \| \theta_j'-\theta_j\|,\quad
\forall \theta_j', \theta_j\in\Theta_j,
\edeq
where the uniformity in $\theta_j$ can be established by finite covering theorem under the compactness of $\Theta_j$.
By virtue of 
\cite[Theorem 1]{klatte1987note}, we have 
\bgeq 
\big|v_j(x_j,x_{-j},\theta_j')- v_j(x_j,x_{-j},\theta_j)\big|
\leq C_1\|\theta_j'-\theta_j\|, \quad \forall \theta_j',\theta_j\in\Theta_j.
\edeq 
The proof is complete.
\hfill $\Box$

With the Proposition \ref{pro:pro1}, we are now ready to address
the existence of BDRNE in the next theorem. 

\begin{theorem}[Existence of an equilibrium]
\label{the:existence}
Under Assumption \ref{ass:Assumption_1}, the following assertions hold.
\begin{itemize}

\item[(i)] $\hat{x}^N_\epsilon$ is a BDRNE that satisfies \eqref{eq:BDRNE} iff $\displaystyle{\hat{x}^N_\epsilon \in \mathop{\arg \max}_{x\in \mathcal{X}}\ \phi^N_\epsilon(x, \hat{x}^N_\epsilon)}$, where $\phi^N_\epsilon(x, \hat{x}^N_\epsilon)$ is defined as  
\bgeqn\label{eq:phi_epsilon}
\phi^N_\epsilon(x, \hat{x}^N_\epsilon):= \sum_{j=1}^n \mathbb{E}_{\theta_j^{N_j}} \left[\min_{Q\in \mathcal{Q}_{\epsilon_j}^{\theta_j}} \mathbb{E}_{Q}\big[u_j(x_j, \hat{x}_{\epsilon_{-j}}^{N_{-j}}, \xi)\big] \right] = \sum_{j=1}^n \vt_j(x_j, \hat{x}_{\epsilon_{-j}}^{N_{-j}}).
\edeqn 
 
\item[(ii)] There exists a BDRNE that satisfies \eqref{eq:BDRNE}. 
 \end{itemize}

\end{theorem}
\noindent
\textbf{Proof.} Under Assumption \ref{ass:Assumption_1}, $\mathbb{E}_{Q}[u_j(x_j, x_{-j}, \xi)]$ is both concave and continuous for every $Q\in \mathcal{Q}_{\epsilon_j}^{\theta_j}$. The concavity of $\mathbb{E}_{Q}[u_j(x_j, x_{-j}, \xi)]$ is preserved by the infimum operation, and the continuity holds under the compactness of $\mathcal{Q}_{\epsilon_j}^{\theta_j}$ (under the topology of weak convergence). Consequently, $\phi^N_\epsilon(\cdot, \hat{x}^N_\epsilon)$ is concave and continuous with respect to $x\in \mathcal{X}$ for any fixed $\hat{x}^N_\epsilon \in \mathcal{X}$. 

\underline{Part (i)}. The reformulation is well known for deterministic Nash equilibrium, see e.g.~\cite{rosen1965existence}. 
Here we provide a proof since the BDRO model \eqref{eq:BDRO} involves maximin operations. The ``if'' part follows from the fact that if $\hat{x}_\epsilon^N$ is not a BDRNE, then there exists $\hat{x}_j \neq \hat{x}_{\epsilon_j}^{N_j}$ such that
$$
\vt_j(\hat{x}_j, \hat{x}_{\epsilon_{-j}}^{N_{-j}}) > 
\vt_j(\hat{x}_{\epsilon_j}^{N_j}, \hat{x}_{\epsilon_{-j}}^{N_{-j}}). 
$$
By defining $\hat{x}:= (\hat{x}_{\epsilon_1}^{N_1}, \cdots, \hat{x}_{\epsilon_{j-1}}^{N_{j-1}}, \hat{x}_j, \hat{x}_{\epsilon_{j+1}}^{N_{j+1}}, \cdots, \hat{x}_{\epsilon_n}^{N_n})$, it follows that $\phi^N_\epsilon(\hat{x}, \hat{x}_\epsilon^N) > \phi^N_\epsilon(\hat{x}_\epsilon^N, \hat{x}_\epsilon^N)$, leading to a contradiction. The ``only if'' part is obvious in that
$$
\vt_j(x_j, \hat{x}_{\epsilon_{-j}}^{N_{-j}}) \leq 
\vt_j(\hat{x}_{\epsilon_j}^{N_j}, \hat{x}_{\epsilon_{-j}}^{N_{-j}}).
$$
Summing over $j$ on both sides of the inequality shows that $\hat{x}_\epsilon^N$ is a BDRNE.   

\underline{Part (ii)}.
The existence of an optimal solution to \eqref{eq:BDRNE} is guaranteed by 
Proposition~\ref{pro:pro1} (ii)-(iii).
To complete the proof, we need to show the existence of $y^\ast\in \mathcal{X}$ such that 
$$
y^\ast \in \mathop{\arg \max}_{x\in \mathcal{X}}\ \phi^N_\epsilon(x, y^\ast)\; \inmat{ or } \;
y^\ast \in \mathop{\arg \min}_{x\in \mathcal{X}} -\phi^N_\epsilon(x, y^\ast). 
$$
Let $\Gamma(y)$ represent the set of optimal solutions to $\min_{x\in \mathcal{X}} -\phi^N_\epsilon(x, y)$ for each fixed $y\in \mathcal{X}$. Then $\Gamma(y) \subset \mathcal{X}$. 
By the convexity of $-\phi^N_\epsilon(\cdot, y)$, $\Gamma(y)$ is a convex set.
Moreover, it can be shown that $\Gamma(y)$ is closed, that is, for $y_k\rightarrow \hat{y}$ and $x_k\in \Gamma(y_k)$ with $x_k\rightarrow \hat{x}$, $\hat{x}\in \Gamma(\hat{y})$. 
Moreover, by \cite[Theorem 4.2.1]{bank1982non}, $\Gamma(\cdot)$ is upper semi-continuous on $\mathcal{X}$. By the fixed point theorem \citep{kakutani1941generalization}, there exists $y^\ast \in \mathcal{X}$ such that $y^\ast \in \Gamma(y^\ast)$.    
\hfill $\Box$

As a special case of BDRNE, we may also reformulate and
simplify the BANE model as a single optimization problem by define
\begin{equation}\label{eq:phi_average}
\phi^N(x, \bar{x}^N) := \sum_{j=1}^n \mathbb{E}_{\theta_j^{N_j}} \left[ \mathbb{E}_{\xi|\theta_j}\big[u_j(x_j, \bar{x}_{-j}^{N_{-j}}, \xi)\big] \right].
\end{equation}
Under Assumption \ref{ass:Assumption_1}, the function \eqref{eq:phi_average} is 
well defined. 
We can establish the existence of BANE analogous to Theorem \ref{the:existence}, we omit the details here.

\subsection{Relation between BDRNE and BANE}

Note that when $\epsilon_j = 0$ in \eqref{eq:ambiguitysetQ_j}, the BDRNE reduces to the BANE. The next proposition quantifies the difference between the two equilibria.


\begin{proposition}\label{pro:StabilityofBDRNE}
Let $\phi^N_\epsilon(x, \hat{x}_\epsilon^N)$ and $\phi^N(x, \bar{x}^N)$ be defined as in \eqref{eq:phi_epsilon} and \eqref{eq:phi_average}, respectively, where $\hat{x}_\epsilon^N$ and $\bar{x}^N$ are the unique BDRNE and BANE that satisfy \eqref{eq:BDRNE} and \eqref{eq:Bayesian-Average}, respectively.
Assume: (a) Assumption~\eqref{ass:Assumption_1} holds, (b) there exist constants $\overline{u}_j$ and $\underline{u}_j$ 
such that $u_j(x_j, x_{-j}, \xi)\in [\underline{u}_j, \overline{u}_j]$.
Then following assertions hold.
\begin{itemize}
    \item [(i)] If $u_j(x_j, x_{-j}, \xi)$ is Lipschitz continuous in $(x_j, x_{-j})$ with modulus $\kappa_j(\xi)$, i.e., 
    \begin{equation}\label{eq:lipschitz}
    \big\| u_j(x'_j, x'_{-j}, \xi) - u_j(x''_j, x''_{-j}, \xi) \big\| \leq \kappa_j(\xi) \| x' - x'' \|,    
    \end{equation}
    then
    \bgeq\label{eq:stabilityanalysis}
    \sup_{x\in \mathcal{X}} \big|\phi^N_\epsilon(x, \hat{x}_\epsilon^N) - \phi^N(x, \bar{x}^N) \big| \leq n \left(\mu  \sqrt{\frac{\epsilon}{2}} + \kappa_N \|\hat{x}_\epsilon^N - \bar{x}^N\|\right),
    \edeq
    where $\mu  := \max_j | \overline{u}_j - \underline{u}_j|$, and
    \bgeqn
    \label{eq:kappa-PJO}
 \kappa_N:=  
\sum_{j=1}^n \mathbb{E}_{\theta_j^{N_j}} \left[ \mathbb{E}_{\xi | \theta_j} \big[ \kappa_j(\xi) \big] \right].
    \edeqn  
    
    \item [(ii)] If, in addition, $\phi^N_\epsilon(\hat{x}_\epsilon^N, \hat{x}_\epsilon^N)$ satisfies a growth condition at $\hat{x}_\epsilon^N$, that is, there exist constants $\gamma > 0$ and $0 < \nu < 1$ such that
    \begin{equation}\label{eq:growthcondition}
    - \phi^N_\epsilon(x, \hat{x}_\epsilon^N) \geq - \phi^N_\epsilon(\hat{x}_\epsilon^N, \hat{x}_\epsilon^N) + \gamma\| x - \hat{x}_\epsilon^N\|^\nu,  \quad \forall x\in \mathcal{X}, 
    \end{equation}
    then 
    $$
    \| \hat{x}_\epsilon^N - \bar{x}^N \| \leq \left(\frac{2\sqrt{2\epsilon}n\mu }{\gamma}\right)^{\frac{1}{\nu}}
    $$
    for sufficient small $\epsilon$ such that $\|\hat{x}_\epsilon^N - \bar{x}^N\|\leq (\frac{\gamma}{2\kappa_N})^{\frac{1}{1 - \nu}}$.
\end{itemize}
\end{proposition}

\noindent
\textbf{Proof.} \underline{Part (i)}. Observe that
\begin{align*}
\sup_{x\in \mathcal{X}} \big |\phi^N_\epsilon(x, \hat{x}_\epsilon^N) - \phi^N(x, \bar{x}^N)\big |
& = 
\sup_{x\in \mathcal{X}} \big |
\phi^N_\epsilon(x, \hat{x}_\epsilon^N) - \phi^N(x, \hat{x}_\epsilon^N)
+ \phi^N(x, \hat{x}_\epsilon^N) -  \phi^N(x, \bar{x}^N) \big| \\
& \leq \sup_{x\in \mathcal{X}} \big |
\phi^N_\epsilon(x, \hat{x}_\epsilon^N) - \phi^N(x, \hat{x}_\epsilon^N) \big| 
+ \sup_{x\in \mathcal{X}} \big | \phi^N(x, \hat{x}_\epsilon^N) -  \phi^N(x, \bar{x}^N) \big| \\
& := R_1 + R_2. 
\end{align*}
We first estimate $R_1$. By definition.
\begin{align*}
R_1 & = \sup_{x \in \X} \left| \sum_{j=1}^n \mathbb{E}_{\theta_j^{N_j}} \left[\min_{Q\in \mathcal{Q}_{\epsilon_j}^{\theta_j}} \mathbb{E}_{Q}\big[u_j(x_j, \hat{x}_{\epsilon_{-j}}^{N_{-j}}, \xi)\big]\right]
- \sum_{j=1}^n \mathbb{E}_{\theta_j^{N_j}} \left[ \mathbb{E}_{\xi | \theta_j} \big[u_j(x_j, \hat{x}_{\epsilon_{-j}}^{N_{-j}}, \xi)\big] \right]
\right| \\
& = \sup_{x \in \X} \left| \sum_{j=1}^n \mathbb{E}_{\theta_j^{N_j}} \left[
\min_{Q\in \mathcal{Q}_{\epsilon_j}^{\theta_j}} \mathbb{E}_{Q}\big[u_j(x_j, \hat{x}_{\epsilon_{-j}}^{N_{-j}}, \xi)\big] -   \mathbb{E}_{\xi | \theta_j} \big[u_j(x_j, \hat{x}_{\epsilon_{-j}}^{N_{-j}}, \xi)\big] \right] \right | \\
& \leq  \sup_{x \in \X}  \sum_{j=1}^n \mathbb{E}_{\theta_j^{N_j}}
\left| \min_{Q\in \mathcal{Q}_{\epsilon_j}^{\theta_j}} \mathbb{E}_{Q}\big[u_j(x_j, \hat{x}_{\epsilon_{-j}}^{N_{-j}}, \xi)\big] -   \mathbb{E}_{\xi | \theta_j} \big[u_j(x_j, \hat{x}_{\epsilon_{-j}}^{N_{-j}}, \xi)\big]  \right | \\
& \leq  \sup_{x \in \X}  \sum_{j=1}^n \mathbb{E}_{\theta_j^{N_j}} \left[
\max_{Q\in \mathcal{Q}_{\epsilon_j}^{\theta_j}} \left| \mathbb{E}_{Q}\big[u_j(x_j, \hat{x}_{\epsilon_{-j}}^{N_{-j}}, \xi)\big] -   \mathbb{E}_{\xi | \theta_j} \big[u_j(x_j, \hat{x}_{\epsilon_{-j}}^{N_{-j}}, \xi)\big]  \right | \right] \\
& \leq  \sum_{j=1}^n \sup_{x_j \in \X_j} \mathbb{E}_{\theta_j^{N_j}} \left[
\max_{Q\in \mathcal{Q}_{\epsilon_j}^{\theta_j}} \left| \mathbb{E}_{Q}\big[u_j(x_j, \hat{x}_{\epsilon_{-j}}^{N_{-j}}, \xi)\big] -   \mathbb{E}_{\xi | \theta_j} \big[u_j(x_j, \hat{x}_{\epsilon_{-j}}^{N_{-j}}, \xi)\big]  \right | \right]\\
& \leq\sum_{j=1}^n  \mathbb{E}_{\theta_j^{N_j}} \left[ \sup_{x_j \in \X_j} 
\max_{Q\in \mathcal{Q}_{\epsilon_j}^{\theta_j}} \left| \mathbb{E}_{Q}\big[u_j(x_j, \hat{x}_{\epsilon_{-j}}^{N_{-j}}, \xi)\big] -   \mathbb{E}_{\xi | \theta_j} \big[u_j(x_j, \hat{x}_{\epsilon_{-j}}^{N_{-j}}, \xi)\big]  \right | \right] \\
& = \sum_{j=1}^n  \mathbb{E}_{\theta_j^{N_j}} \left[
\max_{Q\in \mathcal{Q}_{\epsilon_j}^{\theta_j}} \sup_{x_j \in \X_j} \left| \mathbb{E}_{Q}\big[u_j(x_j, \hat{x}_{\epsilon_{-j}}^{N_{-j}}, \xi)\big] -   \mathbb{E}_{\xi | \theta_j} \big[u_j(x_j, \hat{x}_{\epsilon_{-j}}^{N_{-j}}, \xi)\big]  \right | \right] \\
& = \sum_{j=1}^n  \mathbb{E}_{\theta_j^{N_j}} \left[
\max_{Q\in \mathcal{Q}_{\epsilon_j}^{\theta_j}} \sup_{x_j \in \X_j} \left| \mathbb{E}_{Q}\left[u_j(x_j, \hat{x}_{\epsilon_{-j}}^{N_{-j}}, \xi) - \frac{\overline{u}_j + \underline{u}_j}{2}\right] -   \mathbb{E}_{\xi | \theta_j} \left[u_j(x_j, \hat{x}_{\epsilon_{-j}}^{N_{-j}}, \xi) - \frac{\overline{u}_j + \underline{u}_j}{2} \right]  \right | \right]\\
& \leq \sum_{j=1}^n  \mathbb{E}_{\theta_j^{N_j}} \left[
\max_{Q\in \mathcal{Q}_{\epsilon_j}^{\theta_j}} \sup_{|g| \leq 1} \frac{\overline{u}_j - \underline{u}_j }{2} \cdot \left| \mathbb{E}_{Q}\big[g(\xi)\big] -   \mathbb{E}_{\xi | \theta_j} \big[g(\xi)\big]  \right | \right] \\
& =  \sum_{j=1}^n (\overline{u}_j - \underline{u}_j) \cdot \mathbb{E}_{\theta_j^{N_j}} 
  \left[\max_{Q\in \mathcal{Q}_{\epsilon_j}^{\theta_j}} \dd_{TV}\left(Q, Q_{\theta_j} \right) \right]\\
& \leq  \sum_{j=1}^n  (\overline{u}_j - \underline{u}_j) \cdot  \mathbb{E}_{\theta_j^{N_j}} \left[
\max_{Q\in \mathcal{Q}_{\epsilon_j}^{\theta_j}} \sqrt{ \frac{\dd_{KL}(Q || Q_{\theta_j} )}{2} } \right]\\
& \leq  \sum_{j=1}^n  
 (\overline{u}_j - \underline{u}_j) \cdot \sqrt{\frac{\epsilon}{2}} 
 \leq n \cdot \max_j | \overline{u}_j - \underline{u}_j | \cdot \sqrt{\frac{\epsilon}{2}} = n\mu  \sqrt{\frac{\epsilon}{2}},
\end{align*}
\noindent
where $\mu  := \max_j | \overline{u}_j - \underline{u}_j|$, and $\dd_{TV}$ denotes the total variation distance, and the third last inequality is based on the conclusion in \cite{gibbs2002choosing}.
Next, we estimate $R_2$.
\begin{align*}
R_2 & = \sup_{x\in \mathcal{X}} \left| \sum_{j=1}^n \mathbb{E}_{\theta_j^{N_j}} \left[ \mathbb{E}_{\xi | \theta_j} \big[u_j(x_j, \hat{x}_{\epsilon_{-j}}^{N_{-j}}, \xi)\big] \right] - \sum_{j=1}^n \mathbb{E}_{\theta_j^{N_j}} \left[ \mathbb{E}_{\xi | \theta_j} \big[u_j(x_j, \bar{x}_{-j}^{N_{-j}}, \xi)\big] \right] \right| \\
& = \sup_{x\in \mathcal{X}} \left| \sum_{j=1}^n \mathbb{E}_{\theta_j^{N_j}} \left[ \mathbb{E}_{\xi | \theta_j} \big[u_j(x_j, \hat{x}_{\epsilon_{-j}}^{N_{-j}}, \xi) - u_j(x_j, \bar{x}_{-j}^{N_{-j}}, \xi)\big] \right] \right| \\
& \leq  \sum_{j=1}^n \mathbb{E}_{\theta_j^{N_j}} \left[ \mathbb{E}_{\xi | \theta_j} \big[ \kappa_j(\xi) \|\hat{x}_\epsilon^N - \bar{x}^N\| \big] \right] = \kappa_N \|\hat{x}_\epsilon^N - \bar{x}^N\|,
\end{align*}
where 
the 
inequality is 
due to \eqref{eq:lipschitz}. Thus
\begin{align*}
\sup_{x\in \mathcal{X}} \big |\phi^N_\epsilon(x, \hat{x}_\epsilon^N) - \phi^N(x, \bar{x}^N)\big | \leq  n\mu  \sqrt{\frac{\epsilon}{2}} + \kappa_N \|\hat{x}_\epsilon^N - \bar{x}^N\|.
\end{align*}

\underline{Part (ii)}. We first estimate $\big|\phi^N(\hat{x}_\epsilon^N, \hat{x}_\epsilon^N) - \phi^N(\bar{x}^N, \hat{x}_\epsilon^N)\big|$.
\begin{align*}
& \; \big|\phi^N(\hat{x}_\epsilon^N, \hat{x}_\epsilon^N) - \phi^N(\bar{x}^N, \hat{x}_\epsilon^N)\big| \\
& = \left| \sum_{j=1}^n \mathbb{E}_{\theta_j^{N_j}} \left[ \mathbb{E}_{\xi | \theta_j} \big[u_j(\hat{x}_{\epsilon_j}^{N_j}, \hat{x}_{\epsilon_{-j}}^{N_{-j}}, \xi)\big] \right] - \sum_{j=1}^n \mathbb{E}_{\theta_j^{N_j}} \left[ \mathbb{E}_{\xi | \theta_j} \big[u_j(\bar{x}_j^{N_j}, \hat{x}_{\epsilon_{-j}}^{N_{-j}}, \xi)\big] \right] \right|\\
& = \left| \sum_{j=1}^n \mathbb{E}_{\theta_j^{N_j}} \left[ \mathbb{E}_{\xi | \theta_j} \big[u_j(\hat{x}_{\epsilon_j}^{N_j}, \hat{x}_{\epsilon_{-j}}^{N_{-j}}, \xi) - u_j(\bar{x}_j^{N_j}, \hat{x}_{\epsilon_{-j}}^{N_{-j}}, \xi) \big] \right] \right|\\
& \leq \sum_{j=1}^n \mathbb{E}_{\theta_j^{N_j}} \left[ \mathbb{E}_{\xi | \theta_j} \big[ \kappa_j(\xi) \|\hat{x}_\epsilon^N - \bar{x}^N\| \big] \right] = \kappa_N \|\hat{x}_\epsilon^N - \bar{x}^N\|,
\end{align*}
where the 
last inequality is 
due to \eqref{eq:lipschitz}.
Under the growth condition \eqref{eq:growthcondition}, we can set $x:= \bar{x}^N$ and obtain 
\begin{align*}
 \gamma \| \hat{x}_\epsilon^N - \bar{x}^N \|^\nu 
& \leq \phi^N_\epsilon(\hat{x}_\epsilon^N, \hat{x}_\epsilon^N) - \phi^N_\epsilon(\bar{x}^N, \hat{x}_\epsilon^N) \\
& \leq \big|\phi^N_\epsilon(\hat{x}_\epsilon^N, \hat{x}_\epsilon^N) - \phi^N(\hat{x}_\epsilon^N, \hat{x}_\epsilon^N)\big| + \big|\phi^N(\hat{x}_\epsilon^N, \hat{x}_\epsilon^N) - \phi^N(\bar{x}^N, \hat{x}_\epsilon^N)\big|\\
& \quad + \big|\phi^N(\bar{x}^N, \hat{x}_\epsilon^N) - \phi^N_\epsilon(\bar{x}^N, \hat{x}_\epsilon^N)\big|\\
& \leq 2 \sup_{x \in \X} \left| \phi^N_\epsilon(x, \hat{x}_\epsilon^N) - 
\phi^N(x, \hat{x}_\epsilon^N) \right|  +  \big|\phi^N(\hat{x}_\epsilon^N, \hat{x}_\epsilon^N) - \phi^N(\bar{x}^N, \hat{x}_\epsilon^N)\big|\\
& \leq 2 n \mu \sqrt{\frac{\epsilon}{2}} + \kappa_N \|\hat{x}_\epsilon^N - \bar{x}^N\|. 
\end{align*}
Let $\epsilon$ be sufficiently small such that
$$
\frac{1}{2}\gamma \| \hat{x}_\epsilon^N - \bar{x}^N \|^\nu \geq \kappa_N \| \hat{x}_\epsilon^N - \bar{x}^N \|,
$$
that is, $\|\hat{x}_\epsilon^N - \bar{x}^N\|\leq (\frac{\gamma}{2\kappa_N})^{\frac{1}{1 - \nu}}$. Then we can derive
$\frac{1}{2}\gamma \| \hat{x}_\epsilon^N - \bar{x}^N \|^\nu \leq n\mu  \sqrt{2\epsilon},$
which implies $\|\hat{x}_\epsilon^N - \bar{x}^N \| \leq (\frac{2\sqrt{2\epsilon}n\mu }{\gamma})^{\frac{1}{\nu}}$, and it is obvious that $\hat{x}_\epsilon^N \rightarrow \bar{x}^N$ for $\epsilon \rightarrow 0$.
\hfill $\Box$

\subsection{Convergence of BDRNE and BANE}
We now move on to investigate 
convergence of BDRNE 
as $\min_j N_j \rightarrow \infty$, for $j = 1, \cdots, n$.
Let $Q^\ast$ denote
the true distribution of $\xi$ and 
$q^\ast(\xi)$ 
the corresponding 
density function.
Let
\bgeqn 
\label{eq:Theta^*}
\Theta_j^*:=
\mathop{\arg\min}_{\theta_j \in \Theta_j}
\dd_{KL}\big(Q^*  || Q_{
\theta_j}\big),
\edeqn 
which is a subset 
parameter values in $\Theta_j$ 
such that 
the nominal distribution $Q_{\theta_j^*}$ closest to $Q^\ast$, for $\theta_j^*\in \Theta_j^*$. 
In the case that 
the true distribution $Q^*$ lies in the parametric family 
$\{Q_{\theta_j}\}$,
$\Theta_j^*$ is a singleton $\{\theta_j^*\}$. 
Consequently, the optimal value
of problem (\ref{eq:Theta^*}) is zero 
and 
$f_{\theta_j^*}(\xi)=q^\ast(\xi)$.
In general, $\Theta_j^*$ 
is not a singleton. 
The next example explains this.

\begin{example}
    Consider
the case that $Q^*$ 
is the uniform distribution
over $[0,1]$ and 
$Q_{\theta_j}$ is a family of 
parametric distributions
with density functions
$$
f_{\theta_j}(t) = 
\begin{cases}
\frac{1}{2}, & 
\inmat{for}\; t \in [0, \theta_j-\frac{1}{4}], 
\\
8(t-\theta_j+\frac{1}{4})+\frac{1}{2}, & \inmat{for}\; t \in [\theta_j - \frac{1}{4}, \theta_j], \\
-8(t-\theta_j+\frac{1}{4})+\frac{9}{2}, & \inmat{for}\; t \in [\theta_j, \theta_j + \frac{1}{4}], \\
\frac{1}{2}, & \inmat{for}\; t \in [\theta_j + \frac{1}{4}, 1],
\end{cases}
$$
where $\frac{1}{4} \leq \theta_j \leq \frac{3}{4}$. To facilitate reading, we 
plot 
the density function in Figure \ref{fig:f(t)}.
We can calculate that 
$$
\dd_{KL}\big(Q^\ast || Q_ {\theta_j}\big) = \frac{1}{2} + \ln 2 - \frac{5}{8}\ln 5, \quad \forall \theta_j\in \left[\frac{1}{4}, \frac{3}{4}\right],
$$
which 
implies $\Theta_j^* =\left[\frac{1}{4}, \frac{3}{4}\right]$.
\begin{figure}[H] 
  \centering
\includegraphics[width=0.6\textwidth]{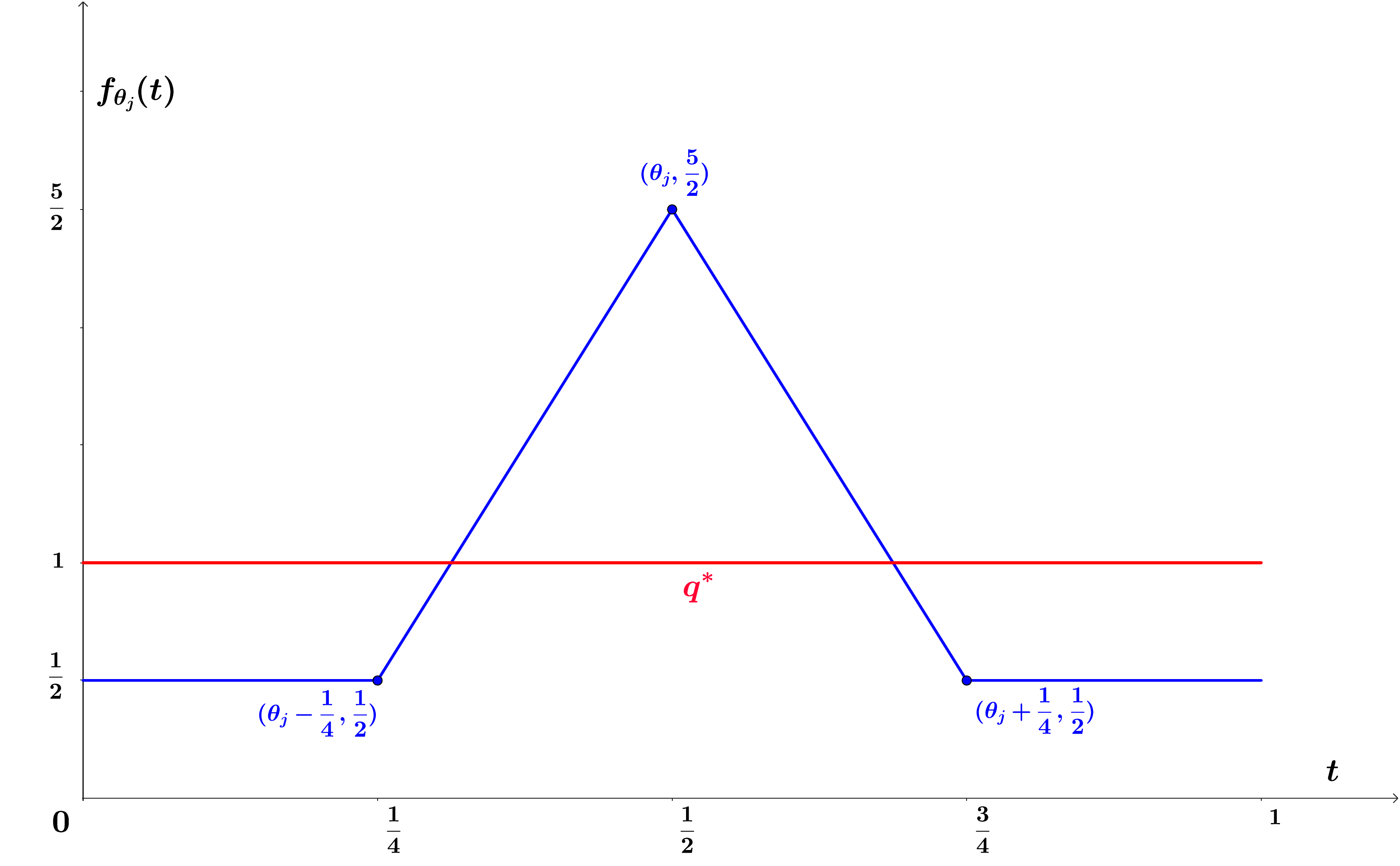}

\captionsetup{font = footnotesize}

  \caption{The blue solid curve 
  represents
  the density function $f_{\theta_j}(t)$, which is equal to $\frac{1}{2}$ when $t\in [0, \theta_j - \frac{1}{4}]\cup [\theta_j+\frac{1}{4}, 1]$, and $f_{\theta_j}(t) = 8(t-\theta_j+\frac{1}{4})+\frac{1}{2}$ when $t\in [\theta_j-\frac{1}{4}, \theta_j]$, $f_{\theta_j}(t) = -8(t-\theta_j+\frac{1}{4})+\frac{9}{2}$ when $t\in [\theta_j, \theta_j+\frac{1}{4}]$. The red solid line 
  represents the true density function $q^\ast$.
  }
\label{fig:f(t)}
\end{figure}
\end{example}
It follows by \cite[Theorem 3.1]{shapiro2023bayesian} and the discussion 
prior to the theorem that the posterior distribution $\mu^{N_j}_{\theta_j}$ concentrates
on a neighborhood of $\epsilon$-optimal solutions of  problem (\ref{eq:Theta^*}) for any fixed $\epsilon$ w.p.1.~(of the distribution of $\Xi\times \Xi\times\cdots$) as $N_j$ goes to infinity.
The next proposition states uniform convergence of $\mathbb{E}_{\theta_j^{N_j}}[v_j(x, \theta_j)]$ to $\mathbb{E}_{\mu_{\theta_j|\Theta_j^*}}[v_j(x_j, \theta_j)]$ as $N_j \rightarrow +\infty$, for $j = 1, \cdots, n$.

\begin{proposition}
\label{pro:LLNBDRNE}
Assume the setting and conditions 
of Proposition~\ref{pro:pro1}, assume also that
$v_j(x, \theta_j)$ is continuous 
in $(x, \theta_j)$.
Let $\mu^{N_j}_{\theta_j}$ denote 
the posterior distribution of $\theta_j$ and
$\mu_{\theta_j|\Theta_j^*}$ denote 
some distribution
of $\theta_j$ over $\Theta_j^*$. 
If $\dd_K(\mu^{N_j}_{\theta_j}, \mu_{\theta_j|\Theta_j^*})\to 0$
 as $N_j\to\infty$, 
then
\begin{equation}\label{eq:v_j-uniform-convg}
\lim_{N_j\rightarrow \infty} \sup_{x\in \mathcal{X}} \left|\mathbb{E}_{\theta_j^{N_j}}\big[v_j(x, \theta_j)\big] - \mathbb{E}_{\mu_{\theta_j|\Theta_j^*}}[v_j(x, \theta_j)]\right| = 0, \; \inmat{w.p.1}, 
\end{equation}
for $j = 1, \cdots, n$, where we write $v_j(x, \theta_j)$ for $v_j(x_j,x_{-j}, \theta_j)$ for simplicity.
\end{proposition}

\noindent
\textbf{Proof.} By Proposition~\ref{pro:pro1} (iv),
$v_j(x, \theta_j)$ is Lipschitz continuous in $\theta_j$ uniformly for all $x\in {\cal X}$.
Thus
\begin{align*}
\left|\mathbb{E}_{\theta_j^{N_j}}\big[v_j(x, \theta_j)\big] - \mathbb{E}_{\mu_{\theta_j|\Theta_j^*}}\big[v_j(x, \theta_j)\big]\right|
& = C_1\left|\mathbb{E}_{\theta_j^{N_j}}\big[v_j(x, \theta_j)/C_1\big] - \mathbb{E}_{\mu_{\theta_j|\Theta_j^*}}\big[v_j(x, \theta_j)/C_1\big]\right|\\
& \leq \sup_{g\in \mathscr{G}}
\left| 
\mathbb{E}_{\theta_j^{N_j}}\big[g\big] - \mathbb{E}_{\mu_{\theta_j|\Theta_j^*}}\big[g\big]\right|
= C_1\dd_K(\mu^{N_j}_{\theta_j}, \mu_{\theta_j|\Theta_j^*}),
\end{align*}
where $\mathscr{G}$ denote the set of Lipschitz continuous function over $\Theta$ with modulus $C_1$.
By taking supremum w.r.t. $x$ on both sides of the inequality above, we obtain (\ref{eq:v_j-uniform-convg}).
\hfill $\Box$

Note that since $\Theta_j$ is a compact set and the Kantorovich distance metricizes the topology of weak convergence, the assumption $\dd_K(\mu^{N_j}_{\theta_j},\mu_{\theta_j|\Theta_j^*})\to 0$
means $\mu^{N_j}_{\theta_j}$ converges weakly to $\mu_{\theta_j|\Theta_j^*}$.
Note 
also 
that in the case when $\Theta_j^*$ is a singleton, 
$\mu_{\theta_j|\Theta_j^*}$
reduces to the Dirac delta distribution of $\theta_j$ at $\theta_j^*$. Consequently, 
we can use \cite[Theorem 3.2]{shapiro2023bayesian}
to establish (\ref{eq:v_j-uniform-convg})
without conditions in Proposition \ref{pro:pro1} (iv).

Proposition \ref{pro:LLNBDRNE} indicates that $\mathbb{E}_{\theta_j^{N_j}}[v_j(x, \theta_j)]$ converges to $\mathbb{E}_{\mu_{\theta_j|\Theta_j^*}}[v_j(x, \theta_j)]$ as the sample size $N_j \rightarrow \infty$ tends to infinity. Consequently, we define
\bgeqn\label{eq:phi_epsilon_infty}
\phi^\infty_\epsilon(x, \hat{x}^\infty_\epsilon):= \sum_{j=1}^n \mathbb{E}_{\mu_{\theta_j|\Theta_j^*}}\left[
\min_{Q\in \mathcal{Q}_{\epsilon_j}^{\theta_j}} 
\mathbb{E}_{Q}\big[u_j(x_j, \hat{x}_{\epsilon_{-j}}^{\infty}, \xi)\big] \right],
\edeqn 
where $\hat{x}^\infty_\epsilon$ represents the BDRNE 
of the following problem:
\begin{equation}
\label{eq:phi_epsilon_infty2}
\hat{x}^\infty_{\epsilon_j} 
\in \mathop{\arg \max}_{x_j\in \mathcal{X}_j}\ \mathbb{E}_{\mu_{\theta_j|\Theta_j^*}}
\left[ \min_{Q\in \mathcal{Q}_{\epsilon_j}^{\theta_j}}\mathbb{E}_{Q}\big[ u_j(x_j, \hat{x}_{\epsilon_{-j}}^{\infty}, \xi) 
\big] \right],\quad \inmat{for}\; j=1,\cdots,n.
\end{equation}
Based on \eqref{eq:phi_epsilon_infty}-\eqref{eq:phi_epsilon_infty2} and Proposition \ref{pro:LLNBDRNE},  
we are ready to investigate the convergence of BDRNE. 

\begin{theorem}[Convergence of BDRNE]\label{pro:ConvergenceofBDRNE}
Let $\{\hat{x}_\epsilon^N\}$ be a sequence of BDRNE
obtained from solving \eqref{eq:BDRNE} and  
$\hat{x}^\infty_\epsilon$ be a cluster point.
Let $\phi^N_\epsilon(x, \hat{x}_\epsilon^N)$ and $\phi^\infty_\epsilon(x, \hat{x}^\infty_\epsilon)$ be defined as in \eqref{eq:phi_epsilon} and \eqref{eq:phi_epsilon_infty}, respectively. 
Assume: (a) all conditions in Proposition \ref{pro:LLNBDRNE} hold,
(b) $u_j(x_j, x_{-j}, \xi)$ satisfies \eqref{eq:lipschitz}. Then following assertions hold.
\begin{itemize}
    \item [(i)] 
    $\hat{x}^\infty_\epsilon$ 
    is a BDRNE which solves problem \eqref{eq:phi_epsilon_infty2}.

    \item [(ii)] If, in addition, $\phi^N_\epsilon(\hat{x}_\epsilon^N, \hat{x}_\epsilon^N)$ satisfies a growth condition at $\hat{x}_\epsilon^N$, that is, there exist constants $\gamma > 0$ and $\nu \in (0, 1)$ such that
    $$
    - \phi^N_\epsilon(x, \hat{x}_\epsilon^N) \geq - \phi^N_\epsilon(\hat{x}_\epsilon^N, \hat{x}_\epsilon^N) + \gamma\| x - \hat{x}_\epsilon^N \|^\nu, \quad \forall x\in \mathcal{X},
    $$
    then 
    \bgeqn 
    \| \hat{x}^N_\epsilon - \hat{x}^\infty_\epsilon \| \leq \left(\frac{ 4 D_K^N}{\gamma}\right)^{\frac{1}{\nu}},
    \edeqn 
    for sufficiently large $N$ 
    with $\|\hat{x}^N_\epsilon - \hat{x}^\infty_\epsilon \|\leq (\frac{\gamma}{2
    \kappa'_N})^{\frac{1}{1 - \nu}}$, where 
\bgeqn
\label{eq:D^N_K}
D^N_K:=\sum_{j = 1}^n C_1\dd_K(\mu^{N_j}_{\theta_j},\mu_{\theta_j|\Theta_j^*})
\quad \inmat{and}\quad
\kappa'_N = 
\sum_{j=1}^n
\mathbb{E}_{\theta_j^{N_j}}\left[ 
\max_{Q\in \mathcal{Q}_{\epsilon_j}^{\theta_j}} 
\mathbb{E}_{Q}\left[
\kappa_j(\xi)\right] \right],
\edeqn
and w.p.1~$\kappa'_N<C_3$ for some positive constant $C_3$.
\end{itemize}
\end{theorem}

In the case that $\Theta_j^*$ is a singleton, 
$\mu_{\theta_j|\Theta_j^*}$ reduces to the Dirac delta distribution of $\theta_j$ at $\theta_j^*$ with probability 1 (\cite[Theorem 3.1]{shapiro2023bayesian}). Consequently 
$\mathbb{E}_{\mu_{\theta_j|\Theta_j^*}} \left[v_j(x, \theta_j)\right] =  v_j(x, \theta_j^*)$.
In that case, we can drop condition (a).


\noindent
\textbf{Proof.} 
\underline{Part (i)}.
Observe that
\begin{align*}
\sup_{x \in \X} \big| \phi^\infty_\epsilon(x, \hat{x}_\epsilon^\infty) - \phi^N_\epsilon(x, \hat{x}_\epsilon^N) \big| 
& = 
\sup_{x\in \mathcal{X}} \big |
\phi^\infty_\epsilon(x, \hat{x}_\epsilon^\infty) - \phi^N_\epsilon(x, \hat{x}_\epsilon^\infty) + 
\phi^N_\epsilon(x, \hat{x}_\epsilon^\infty) - \phi^N_\epsilon(x, \hat{x}_\epsilon^N) \big| \\
& \leq \sup_{x\in \mathcal{X}} \big |
\phi^\infty_\epsilon(x, \hat{x}_\epsilon^\infty) - \phi^N_\epsilon(x, \hat{x}_\epsilon^\infty) \big| 
+ \sup_{x\in \mathcal{X}} \big | \phi^N_\epsilon(x, \hat{x}_\epsilon^\infty) - \phi^N_\epsilon(x, \hat{x}_\epsilon^N) \big| \\
& := R_1 + R_2. 
\end{align*}
We first estimate $R_1$. By definition,
\begin{align*}
R_1 & = \sup_{x\in \mathcal{X}} \left |
\sum_{j=1}^n \mathbb{E}_{\mu_{\theta_j|\Theta_j^*}}\left[ v_j(x_j, \hat{x}_{\epsilon_{-j}}^{\infty}, \theta_j)\right] - \sum_{j=1}^n \mathbb{E}_{\theta_j^{N_j}} \left[v_j(x_j, \hat{x}_{\epsilon_{-j}}^{\infty}, \theta_j) \right] \right|\\
& = \sup_{x\in \mathcal{X}} \left |
\sum_{j=1}^n \left( \mathbb{E}_{\mu_{\theta_j|\Theta_j^*}}\left[ v_j(x_j, \hat{x}_{\epsilon_{-j}}^{\infty}, \theta_j) \right] - \mathbb{E}_{\theta_j^{N_j}}  \left[v_j(x_j, \hat{x}_{\epsilon_{-j}}^{\infty}, \theta_j) \right] \right) \right| \\
& \leq 
\sup_{x\in \mathcal{X}}  \sum_{j=1}^n \left| \mathbb{E}_{\mu_{\theta_j|\Theta_j^*}}\left[ v_j(x_j, \hat{x}_{\epsilon_{-j}}^{\infty}, \theta_j) \right] - \mathbb{E}_{\theta_j^{N_j}}  \left[v_j(x_j, \hat{x}_{\epsilon_{-j}}^{\infty}, \theta_j) \right] \right| \\
&\leq 
\sum_{j=1}^n \sup_{x_j\in \mathcal{X}_j}   \left| \mathbb{E}_{\mu_{\theta_j|\Theta_j^*}}\left[ v_j(x_j, \hat{x}_{\epsilon_{-j}}^{\infty}, \theta_j) \right] - \mathbb{E}_{\theta_j^{N_j}}  \left[v_j(x_j, \hat{x}_{\epsilon_{-j}}^{\infty}, \theta_j) \right] \right|
\\
& \leq \sum_{j = 1}^n C_1\dd_K(\mu^{N_j}_{\theta_j},\mu_{\theta_j|\Theta_j^*}) = D_K^N,
\end{align*}
where the last inequality is based on Proposition~\ref{pro:LLNBDRNE}.
Next, we estimate $R_2$.
\begin{align*}
R_2 & = \sup_{x\in \mathcal{X}} \left| \sum_{j=1}^n \mathbb{E}_{\theta_j^{N_j}} \left[\min_{Q\in \mathcal{Q}_{\epsilon_j}^{\theta_j}} \mathbb{E}_{Q}\big[u_j(x_j, \hat{x}_{\epsilon_{-j}}^{\infty}, \xi)\big] \right] 
- \sum_{j=1}^n \mathbb{E}_{\theta_j^{N_j}} \left[\min_{Q\in \mathcal{Q}_{\epsilon_j}^{\theta_j}} \mathbb{E}_{Q}\big[u_j(x_j, \hat{x}_{\epsilon_{-j}}^{N}, \xi)\big] \right] \right| \\
& \leq \sup_{x \in \X} \sum_{j=1}^n \mathbb{E}_{\theta_j^{N_j}} \left[ \max_{Q\in \mathcal{Q}_{\epsilon_j}^{\theta_j}} \left | \mathbb{E}_{Q}\big[u_j(x_j, \hat{x}_{\epsilon_{-j}}^{\infty}, \xi)
- u_j(x_j, \hat{x}_{\epsilon_{-j}}^{N}, \xi)\big] \right| \right] \\
& \leq \sum_{j=1}^n
\mathbb{E}_{\theta_j^{N_j}}\left[ 
\max_{Q\in \mathcal{Q}_{\epsilon_j}^{\theta_j}} 
\mathbb{E}_{Q}\left[
\kappa_j(\xi)\right] \right]
\|\hat{x}^\infty_\epsilon - \hat{x}^N_\epsilon \| =  \kappa'_N \|\hat{x}^\infty_\epsilon - \hat{x}^N_\epsilon\|,
\end{align*}
where the 
last inequality is 
due to \eqref{eq:lipschitz}
and $\kappa_N'$ is bounded 
by a positive constant w.p.1.
Thus
\begin{equation}\label{eq:upperbound}
\sup_{x \in \X} \big| \phi^N_\epsilon(x, \hat{x}_\epsilon^N) - \phi^\infty_\epsilon(x, \hat{x}_\epsilon^\infty) \big| \leq D_K^N + \kappa'_N \|\hat{x}^\infty_\epsilon - \hat{x}^N_\epsilon\| \rightarrow 0, \quad\inmat{as}\ \min_j N_j \rightarrow \infty.    
\end{equation}
Moreover, for any $x^N\in \mathcal{X}$, we have
$\phi^N_\epsilon(x^N, \hat{x}_\epsilon^N) \leq
\phi^N_\epsilon(\hat{x}_\epsilon^N, \hat{x}_\epsilon^N)$. By taking a subsequence if necessary, we assume for the simplicity of presentation that $x^N \rightarrow x$. Further, by \eqref{eq:upperbound} and the continuity of $\phi$, we obtain
$\phi^\infty_\epsilon(x, \hat{x}_\epsilon^\infty) \leq
\phi^\infty_\epsilon(\hat{x}_\epsilon^\infty, \hat{x}_\epsilon^\infty)$. Thus, we arrive at $\hat{x}_\epsilon^\infty \in \mathop{\arg \max}_{x \in \X}\ \phi^\infty_\epsilon(x, \hat{x}_\epsilon^\infty)$,
which implies $\hat{x}^\infty_\epsilon$ is an equilibrium of problem \eqref{eq:phi_epsilon_infty}.

\underline{Part (ii)}. 
We first estimate $\big|\phi^{\infty}_\epsilon(\hat{x}_\epsilon^N, \hat{x}_\epsilon^N) - \phi_\epsilon^\infty (\hat{x}_\epsilon^\infty, \hat{x}_\epsilon^N)\big|$.
\begin{align*}
& \; \big|\phi^{\infty}_\epsilon(\hat{x}_\epsilon^N, \hat{x}_\epsilon^N) - \phi_\epsilon^\infty (\hat{x}_\epsilon^\infty, \hat{x}_\epsilon^N)\big| \\
& = \left| \sum_{j=1}^n \mathbb{E}_{\mu_{\theta_j|\Theta_j^*}}\left[ \min_{Q\in \mathcal{Q}_{\epsilon_j}^{\theta_j}} \mathbb{E}_{Q}\big[u_j(\hat{x}_{\epsilon_j}^{N_j}, \hat{x}_{\epsilon_{-j}}^{N_{-j}}, \xi)\big] \right]
- \sum_{j=1}^n \mathbb{E}_{\mu_{\theta_j|\Theta_j^*}}\left[ \min_{Q\in \mathcal{Q}_{\epsilon_j}^{\theta_j}} \mathbb{E}_{Q}\big[u_j(\hat{x}_{\epsilon_j}^\infty, \hat{x}_{\epsilon_{-j}}^{N_{-j}}, \xi)\big] \right] \right|\\
& \leq \sum_{j=1}^n \mathbb{E}_{\mu_{\theta_j|\Theta_j^*}}\left[ \max_{Q \in Q_{\theta_j}} 
\left| \mathbb{E}_Q \big[u_j(\hat{x}_{\epsilon_j}^{N_j}, \hat{x}_{\epsilon_{-j}}^{N_{-j}}, \xi)
- u_j(\hat{x}_{\epsilon_j}^\infty, \hat{x}_{\epsilon_{-j}}^{N_{-j}}, \xi)\big]  \right| \right] \\
& \leq \sum_{j=1}^n
\mathbb{E}_{\mu_{\theta_j|\Theta_j^*}} \left[ 
\max_{Q\in \mathcal{Q}_{\theta_j}} 
\mathbb{E}_{Q}\left[
\kappa_j(\xi)\right] \right]
\|\hat{x}^N_\epsilon - \hat{x}^\infty_\epsilon \| = \kappa'_N \|\hat{x}^N_\epsilon - \hat{x}^\infty_\epsilon\|,
\end{align*}
where the second last inequality is according to \eqref{eq:lipschitz}.
Under the growth condition \eqref{eq:growthcondition}, we can set $x:=  \hat{x}^\infty_\epsilon$ and obtain 
\begin{align*}
\gamma \| \hat{x}^N_\epsilon - \hat{x}^\infty_\epsilon \|^\nu 
& \leq \phi_\epsilon^N(\hat{x}_\epsilon^N, \hat{x}_\epsilon^N) - \phi_\epsilon^N(\hat{x}_\epsilon^\infty, \hat{x}_\epsilon^N) \\
& \leq  \big|\phi_\epsilon^N(\hat{x}_\epsilon^N, \hat{x}_\epsilon^N) - \phi^{\infty}_\epsilon(\hat{x}_\epsilon^N, \hat{x}_\epsilon^N) \big| + \big|\phi^{\infty}_\epsilon(\hat{x}_\epsilon^N, \hat{x}_\epsilon^N) - \phi_\epsilon^\infty (\hat{x}_\epsilon^\infty, \hat{x}_\epsilon^N)\big| \\
& \quad + \big|\phi_\epsilon^{\infty}(\hat{x}_\epsilon^\infty, \hat{x}^N) - \phi_\epsilon^N(\hat{x}_\epsilon^\infty, \hat{x}_\epsilon^N)\big|\\
& \leq 2 \sup_{x \in \X} \left| \phi^N_\epsilon(x, \hat{x}_\epsilon^N)  - 
\phi^\infty_\epsilon(x, \hat{x}_\epsilon^N) \right|  +  \big|\phi^{\infty}_\epsilon(\hat{x}_\epsilon^N, \hat{x}_\epsilon^N) - \phi_\epsilon^\infty (\hat{x}_\epsilon^\infty, \hat{x}_\epsilon^N)\big| \\
& \leq 2D_K^N + \kappa'_N \|\hat{x}^N_\epsilon - \hat{x}^\infty_\epsilon\|. 
\end{align*}
Let $\min_j N_j$ be sufficiently large such that
$$\frac{1}{2}\gamma \| \hat{x}^N_\epsilon - \hat{x}^\infty_\epsilon \|^\nu \geq \kappa'_N \| \hat{x}^N_\epsilon - \hat{x}^\infty_\epsilon \|,$$
that is, $\|\hat{x}^N_\epsilon - \hat{x}^\infty_\epsilon\|\leq (\frac{\gamma}{2\kappa'_N})^{\frac{1}{1 - \nu}}$. Then we can derive
$\frac{1}{2}\gamma \| \hat{x}^N_\epsilon - \hat{x}^\infty_\epsilon \|^\nu \leq 2 D_K^N$,
which implies 
$\| \hat{x}^N_\epsilon - \hat{x}^\infty_\epsilon \| \leq \left(\frac{4 D_K^N}{\gamma}\right)^{\frac{1}{\nu}}$.
Moreover, if $\dd_K(\mu^{N_j}_{\theta_j}, \mu_{\theta_j|\Theta_j^*})\to 0$
as $N_j\to\infty$. Then $\lim_{N\rightarrow \infty} \| \hat{x}^N_\epsilon - \hat{x}^\infty_\epsilon \| = 0$, that is, $\hat{x}^N_\epsilon$ converges to $\hat{x}^\infty_\epsilon$ with probability 1 as $\min_j N_j \rightarrow \infty$.
\hfill $\Box$


As 
discussed in the previous section, when $\epsilon_j = 0$ in \eqref{eq:ambiguitysetQ_j}, BDRNE reduces to BANE. 
Moreover, Proposition \ref{pro:LLNBDRNE} indicates that $\mathbb{E}_{\theta_j^{N_j}}\big[v_j(x, \theta_j)\big]$ converges to $\mathbb{E}_{\mu_{\theta_j|\Theta_j^*}}[v_j(x, \theta_j)]$ as the sample size $N_j$ tends to infinity.
Consequently, we define
\begin{equation}
\label{eq:phiAast}
\phi^{\infty}(x, x^\ast) := \sum_{j=1}^n \mathbb{E}_{\mu_{\theta_j|\Theta_j^*}}\left[ \mathbb{E}_{\xi|\theta_j}
\big[u_j(x_j, x^\ast_{-j}, \xi)\big] \right], 
\end{equation}
where $x^\ast$ 
denotes the BANE/stochastic NE that satisfies 
\begin{equation}
\label{eq:phiAast2}
x_j^{\ast} \in \mathop{\arg \max}_{x_j\in \mathcal{X}_j}\ 
\mathbb{E}_{\mu_{\theta_j|\Theta_j^*}} 
\left[
\mathbb{E}_{\xi| \theta_j}
\big[ u_j(x_j, x^{\ast}_{-j}, \xi) \big]
\right], 
\quad \inmat{for}\; j=1,\cdots,n.  
\end{equation}
Based on 
\eqref{eq:phiAast}-\eqref{eq:phiAast2}, 
we can state 
convergence of BANE as $\min_j N_j \rightarrow \infty$.

\begin{corollary}[Convergence of BANE]
Let $\{\bar{x}^N\}$ be a sequence of BANE
obtained from solving \eqref{eq:Bayesian-Average} and  
$x^\ast$ be a cluster point.
Let $\phi^N(x, \bar{x}^N)$ and $\phi^\infty(x, x^\ast)$ be defined as in \eqref{eq:phi_average} and \eqref{eq:phiAast}, respectively. Assume: (a) all conditions in Proposition \ref{pro:LLNBDRNE} hold, (b) $u_j(x_j, x_{-j}, \xi)$ satisfies \eqref{eq:lipschitz}. Then following assertions hold.
\begin{itemize}
\item [(i)] $x^\ast$ is a NE which solves problem \eqref{eq:phiAast2}. 
    
    \item [(ii)] Let $D^N_{K}$ and $\kappa_N$
    be defined as in (\ref{eq:D^N_K}) and 
    (\ref{eq:kappa-PJO}) respectively.
    If, in addition, $\phi^N(\bar{x}^N, \bar{x}^N)$ satisfies a growth condition at $\bar{x}^N$, that is, there exist constants $\gamma > 0$ and $0 < \nu < 1$ such that
    $$
    - \phi^N(x, \bar{x}^N) \geq - \phi^N(\bar{x}^N, \bar{x}^N) + \gamma\| x - \bar{x}^N \|^\nu, \quad \forall x\in \mathcal{X},
    $$
    then 
    $$
    \| \bar{x}^N - x^\ast \| \leq \left(\frac{ 4 D^N_{K}}{\gamma}\right)^{\frac{1}{\nu}},
    $$
    for sufficient large $\min_j N_j$ such that $\|\bar{x}^N - x^\ast\|\leq (\frac{\gamma}{2\kappa_N})^{\frac{1}{1 - \nu}}$.
    Moreover, $\bar{x}^N$ converges to $x^\ast$ with probability 1 as $\min_j N_j \rightarrow \infty$. 
\end{itemize}
\end{corollary}

The conclusion follows from 
Theorem \ref{pro:ConvergenceofBDRNE} by setting $\epsilon_j=0$. The only difference is the intermediate 
quantity $\kappa_N$ from $\kappa'_N$ where the latter depends on $\epsilon_j$ whereas the former doesn't.

\subsection{Algorithm for solving BDRNE}

We now turn to discuss how to solve 
the BDRNE problem (\ref{eq:BDRNE}). Since
both the posterior distribution $\mu_{\theta_j}^{N_j}$ and 
the distribution $Q_{\theta_j}$ are continuously distributed in general,
it is difficult to obtain a closed form of the objective function. This prompts us to consider sample average approximation of the posterior distribution and the nominal distribution of the ambiguity set $\mathcal{Q}_{\epsilon_j}^{\theta_j}$.
Due to the specific structure of the ambiguity set, we may derive dual formulation of the inner minimization problem and then formulate the objective function of each player 
as a single minimization problem.


By 
\citep{ben1987penalty, shapiro2017distributionally}, we can reformulate the inner 
minimization problem of \eqref{eq:BDRO} as:
\begin{equation}\label{eq:duality}
\begin{aligned}
\min_{Q\in \mathcal{Q}_{\epsilon_j}^{\theta_j}} \mathbb{E}_{Q}\big[u_j(x_j, x_{-j}, \xi)\big]& = 
- \max_{Q\in \mathcal{Q}_{\epsilon_j}^{\theta_j}} \mathbb{E}_{Q}\big[-u_j(x_j, x_{-j}, \xi)\big] \\
& = -\min_{\lambda > 0} \left\{ \lambda \epsilon_j + \lambda \ln \mathbb{E}_{
\xi | \theta_j }\left[ \exp\left( \frac{-u_j(x_j, x_{-j}, \xi)}{\lambda} \right) \right] \right\}.
\end{aligned}    
\end{equation}
Assume that for all $x \in \X,\ j = 1, \cdots n$,
$$
\mathbb{E}_{
\xi | \theta_j }
\left[ 
\exp
\left( 
\frac{-u_j(x, \xi)}{t} \right) 
\right] < +\infty, \quad
\forall t \in \R, 
\theta_j \in \Theta_j.
$$
Then problem \eqref{eq:duality} is strictly convex in $\lambda >0$.
Consequently, the BDRO problem \eqref{eq:BDRO} can be 
recast as:
\begin{equation}\label{eq:duality_1}
-\min_{x_j\in \mathcal{X}_j} \mathbb{E}_{\theta_j^{N_j}} \left[ \min_{\lambda > 0} \left\{ \lambda \epsilon_j + \lambda \ln \mathbb{E}_{\xi | \theta_j }\left[ \exp\left( \frac{-u_j(x_j, x_{-j}, \xi)}{\lambda} \right) \right] \right\} \right].
\end{equation}
Problem \eqref{eq:duality_1} 
can be viewed as a two-stage stochastic program where
$x_j$ is a first-stage decision variable whereas 
$\lambda$ is a second stage decision variable
to be decided after observation of $\theta_j$.
The sample average approximation of the problem can be formulated as
\begin{equation}
\label{eq:SAA}
-\min_{x_j\in \mathcal{X}_j, \lambda_i > 0} \left\{ \frac{1}{N_{\theta_j}} \sum_{i=1}^{N_{\theta_j}} \left( \lambda_i \epsilon_j + \lambda_i \ln \left(\frac{1}{N_\xi}
\sum_{ \hat{\xi}_{\theta_j^i}^k \in \mathcal{D}_{\theta_j^i}} \exp\left( \frac{-u_j
(x_j, x_{-j}, \hat{\xi}_{\theta_j^i}^k)}{\lambda_i}\right) \right) \right) \right\}, 
\end{equation}
where $\theta_j^i, i=1,\cdots, N_{\theta_j}$ 
are i.i.d.~sampling of $\theta_j$ according to the posterior distribution $\mu_{\theta_j}^{N_j}$ and  $\mathcal{D}_{\theta_j^i}=\{\hat{\xi}_{\theta_j^i}^k\}_{k=1}^{N_\xi}$ is i.i.d.~sampling of $\xi$
according to the distribution of $Q_{\theta_j^i}$.
Convergence of problem \eqref{eq:SAA}
to \eqref{eq:duality_1} 
in terms of the optimal values 
can be established 
similar to 
that of 
\cite[Theorem 5.1]{delage2022shortfall}. We omit the details. 

We use 
Gauss-Seidel-type 
iterative 
scheme \citep{dafermos1983iterative, facchinei2010generalized}
to 
find an approximate BDRNE.
Algorithm \ref{alg:BDRNE} describes the computational procedures.

\begin{algorithm}[H] 
\caption{
Gauss-Seidel-type method for solving 
BDRNE approximately
}
\begin{algorithmic}\label{alg:BDRNE}
\STATE \textbf{input}: 
Generate i.i.d.~samples
$\{\theta_j^i\}_{i=1}^{N_{\theta_j}}$ of $\theta_j$ according to the posterior distribution $\mu_{\theta_j}^{N_j}$ 
and i.i.d.~samples
$\mathcal{D}_{\theta_j^i}=\{\hat{\xi}_{\theta_j^i}^k\}_{k=1}^{N_\xi}$ of $\xi$
according to the distribution of $Q_{\theta_j^i}$ for $i = 1, \cdots, N_{\theta_j}$.
Initial point $x^0:= (x_1^0, \cdots, x_n^0)$ and tolerance $\tau$. 


\STATE \textbf{while $\|x^{k} - x^{k-1}\|^2 > \tau$:} 

\STATE \quad let $k:= k+1$, do the following repeat iteration; 

\STATE \quad \textbf{for $j = 1,\cdots, n$:} 
\begin{adjustwidth}{2em}{0em}
with the fixed $\hat{x}_{-j}^{k} := (x_1^{k}, \cdots, x_{j-1}^{k}, x_{j+1}^{k-1}, \cdots, x_n^{k-1})$ and the ambiguity set size $\epsilon_j$, solve \eqref{eq:SAA}
and denote the solution as $x^{k}_j$, then define $\hat{x}^k := (x_1^k, \cdots, x_{j-1}^k, x_j^k, x_{j+1}^{k-1}, \cdots, x_n^{k-1})$.
\end{adjustwidth}

\STATE \quad\textbf{end} 
 
\STATE \quad
denote the optimal solution of the $k$-th iteration as $x^k:= (x_1^k, \cdots, x_n^k)$.
\STATE \textbf{end}

\STATE \textbf{output}: The approximate optimal solution of BDRNE is $x^\ast:= x^k$.
\end{algorithmic}
\end{algorithm}

Under 
the
convexity 
of the objective function, 
every limit point of the sequence generated by the Gauss-Seidel-type method is an approximation solution of the equilibrium (see e.g.~\citep{facchinei2010generalized}).

\section{BDRNE in Price Competition under MNL demand}\label{se:MNL}
In this section, we apply the BDRNE  model 
to a price competition problem.
Consider an oligopoly market with $n$ substitutable differentiated products, where each product is offered by a distinct firm. In this market, customers' utility from choosing product $j$ is determined by the following equation:
\begin{equation} \label{eq:U_j-MNL}
U_j = \beta^{T} x_j - \alpha p_j + \varepsilon_j, \quad \inmat{for} \;  j=1, \cdots, n,
\end{equation}
where $\beta\in \R^m$ is a vector representing the customer's preference or taste for the observed product characteristic $x_j\in \R^m_+$, $\alpha$ is a positive scalar representing the customer's sensitivity to price $p_j$, and $\varepsilon_j: (\Omega, {\cal F}, \mathbb{P})\to \R$ 
is a random variable representing
idiosyncratic product-specific random demand shock and/or 
customer's utility of unobserved characteristics of product $j$.
In this setup, we define $j = 0$ as the outside option, representing the case that customers choose not to purchase any of the products in this market,
and we set the customer's utility from choosing the outside option as $U_0 = \varepsilon_0$. Assuming that $\varepsilon_j$, $j=0, 1, \cdots,n$ are i.i.d. 
and follow the standard Gumbel distribution $\text{Gumbel}(0,1)$, with location parameter 0 and scale parameter 1,
for fixed 
$\xi:= (\beta, \alpha)$, we can derive
the purchase probability (or called market share) of product $j$ as
\begin{equation}\label{eq:q_j(p, xi)}
q_j(p_j, p_{-j}, \xi) := \mathbb{P}(U_j > U_k, \forall k \neq j) = 
\frac{e^{\beta^T x_j  - \alpha p_j}}{\sum_{k=1}^n e^{\beta^T x_k-\alpha p_k}},   
\end{equation}
for $j=1,\cdots,n$, where $p_{-j}:= \left(p_1, \cdots, p_{j-1}, p_{j+1}, \cdots, p_n \right)$
denotes the vector of prices of all firms except product $j$, and the no-purchase probability is
$$ q_0(p, \xi):
= \mathbb{P}(U_0 > U_j, j = 1, \cdots, n) = \frac{1}{\sum_{k=1}^n e^{\beta^T x_k-\alpha p_k}}.$$
This approach is known as the multinomial logit (MNL) model \citep{mcfadden1973conditional, train2009discrete} and is popular due to its computational traceability. 
Based on \eqref{eq:q_j(p, xi)}, we define
$$u_j(p_j, p_{-j}, c_j, \xi):= (p_j - c_j)q_j(p_j, p_{-j}, \xi) $$
as the normalized profit function of firm $j$, where $c_j$ represents product $j$'s marginal cost.

In practice, customer's preferences for the characteristics and prices of products 
may vary across different customers/sectors.
In such cases, $\xi$ is not constant, rather 
it may be 
a random 
vector 
which 
captures 
customers' random preferences. 
This 
kind of modeling approach
may be traced back to 
earlier work by  Berry et al.~\cite{berry1995automobile}, where the authors propose a random coefficient logit demand model and 
provide a 
model to estimate heterogeneous consumers’ preferences over observed product characteristics in a differentiated product setting. 
Here we consider the case that 
the true probability distribution of 
$\xi$ 
can be approximated by 
a family of 
$\{Q_\theta\}$
parameterized by $\theta$, where
$\theta$ 
follows a prior distribution $\mu_\theta$.
Then 
the expected purchase probability 
of product $j$ is
\bgeq
\label{eq:MNL-mrk-shr}
    q_j(p_j, p_{-j}) :=
\int_{\Theta_j}    \int_{\left\{\xi, \varepsilon | U_j > U_k, \forall k \neq j \right\}} G( d \varepsilon)  Q_{\theta_j}(d \xi) \mu_{\theta_j}(d \theta_j)
    = \mathbb{E}_{\mu_{\theta_j}}\big[ \mathbb{E}_{Q_{\theta_j}}\big[q_j(p_j, p_{-j}, \xi)\big] \big].
\edeq
We may formulate firm $j$'s maximum expected profit pricing problem as 
\bgeq
\label{eq:price-j}
\max_{p_j} \; \mathbb{E}_{\mu_{\theta_j}}\left[ \mathbb{E}_{Q_{\theta_j}}\big [u_j(p_j, p_{-j}, c_j, \xi)\big]\right],  \quad \inmat{for} \;  j=1, \cdots, n.
\edeq
In some data-driven problems, 
the true 
distribution of 
$\xi$
is often unknown and needs to 
be 
estimated
with available empirical data. 
This motivates us to apply the proposed BDRNE model in this setting.
Moreover, we assume that each player sets its price based on its own marginal cost $c_j$, and $c_j \in \mathcal{C}_j$, then $c:= \left(c_1, \cdots, c_n\right)$ is a random vector with support set $\mathcal{C}:=\mathcal{C}_1\times \cdots \times \mathcal{C}_n$. Subsequently, we propose the BDRNE in price competition under MNL demand (BDRNE-MNL) defined as follows.
\begin{definition}[BDRNE-MNL]\label{def:BDRNE-MNL} 
Let $c \in \mathcal{C}$ be a fixed parameter.
The BDRNE-MNL is an $n$-tuple $p^\ast(c):= (p_1^\ast(c_1), \cdots, p_n^\ast(c_n) )\in \mathcal{P}(c):= (\mathcal{P}_1(c_1), \cdots, \mathcal{P}_n(c_n))$ such that
\begin{equation}\label{eq:BDRNE-MNL}
p_j^\ast(c_j) \in \mathop{\arg \max}_{p_j\in \mathcal{P}_j(c_j)}\ \mathbb{E}_{\theta_j^{N_j}} \left[ \min_{Q\in \mathcal{Q}_{\epsilon_j}^{\theta_j}} \mathbb{E}_{Q} \big[u_j(p_j, p^\ast_{-j}(c_{-j}), c_j, \xi)\big] \right],
\end{equation}
where $\mathcal{P}_j(c_j)$ denotes firm $j$'s action space, and $\mathcal{Q}_{\epsilon_j}^{\theta_j}$ is defined as in \eqref{eq:ambiguitysetQ_j}, for $j = 1, \cdots, n$.
\end{definition}

In the setup of \eqref{eq:BDRNE-MNL}, each player determines its best response strategy under the Nash conjecture, assuming that its rivals' best response strategy $p^\ast_{-j}(c_{-j})$ is fixed. A noteworthy distinction between the BDRNE-MNL model and existing Nash equilibrium models is the explicit indication of dependence of $p_j$ on $c_j$, which we refer to as a parameterized Nash equilibrium \citep{feinstein2022continuity}. 
Moreover, the action space ${\cal P}_j$ is contingent on $c_j$, indicating that a firm's range of feasible actions is determined by its type parameter.

\subsection{Existence of BDRNE-MNL}

In the following discussions, we will demonstrate the existence of BDRNE-MNL.
According to Theorem \eqref{the:existence}, BDRNE-MNL exists under the  Assumption \ref{ass:Assumption_1}. We will verify that $u_j(p_j, p_{-j}(c_{-j}), c_j, \xi)$ in \eqref{eq:BDRNE-MNL} satisfies the requirements of Assumption \ref{ass:Assumption_1}. Firstly, the concavity of $u_j(p_j, p_{-j}(c_{-j}), c_j, \xi)$ is guaranteed within a specific interval, as detailed in the proposition below. 
\begin{proposition}[Concavity of $u_j(p_j, p_{-j}(c_{-j}), c_j, \xi)$]\label{pro:concavity_u_j} For each fixed $\xi$, the function $u_j(p_j, p_{-j}(c_{-j}), c_j, \xi)$ is concave in $p_j$ within the range $[c_j, \bar{p}_j(c_j, p_{-j}(c_{-j}), \xi)]$, where $\bar{p}_j(c_j, p_{-j}(c_{-j}), \xi)$ is the unique solution of
\begin{equation}\label{eq:phi_j-p_j-p_{-j}} 
\phi_{j}\big(p_j, p_{-j}(c_{-j}), c_j, \xi \big) := p_j - c_j - \frac{1}{\alpha\big(\frac{1}{2} - q_j(p_j,p_{-j}(c_{-j}), \xi)\big)} = 0 
\end{equation}
over the interval $[c_j, +\infty)$
with $q_j(\bar{p}_j(c_j, p_{-j}(c_{-j}), \xi), p_{-j}(c_{-j}), \xi) < \frac{1}{2}$,
for $j = 1, \cdots, n$.
\end{proposition}
\noindent
\textbf{Proof.}
As stated in Liu et al.~\cite[Lemma 3.2 (iii)]{liusunxu2024Bayesian}, equation \eqref{eq:phi_j-p_j-p_{-j}}
has a unique solution $\bar{p}_j(c_j, p_{-j}(c_{-j}), \xi)$ over the interval $[c_j, +\infty)$. Likewise, by Liu et al.~\cite[Proposition 1 (i)]{liusunxu2024Bayesian}, we demonstrate that function $u_j(p_j, p_{-j}(c_{-j}), c_j, \xi)$ is concave in $p_j$ over the interval $[c_j, \bar{p}_j(c_j, p_{-j}(c_{-j}), \xi)]$.
\hfill $\Box$

In Assumption \ref{ass:Assumption_1} (a), we require $u_j(\cdot, x_{-j}, \xi)$ to be concave over the interval $\mathcal{X}_j$. To ensure this, we define an interval $\mathcal{P}_j(c_j) := [c_j, \bar{p}_j(c_j, p_{-j}(c_{-j}), \hat{\xi})]$ such that $u_j(p_j, p_{-j}(c_{-j}), c_j, \xi)$ is concave in $p_j$ over $\mathcal{P}_j(c_j)$ uniformly for all $\xi$. The following lemma characterizes the properties of $\bar{p}_j(c_j, p_{-j}(c_{-j}), \xi)$.

\begin{lemma}[Property of $\bar{p}_j(c_j, p_{-j}(c_{-j}), \xi)$] \label{lem:barpj}
For each $j = 1, \cdots, n$, (a) $\frac{\partial \bar{p}_j}{\partial \beta} > 0$; (b) $\frac{\partial \bar{p}_j}{\partial \alpha} < 0$.
\end{lemma}
\noindent
\textbf{Proof.}
Using the fact that $\frac{\partial \phi_j}{\partial p_j} > 0$, $\frac{\partial \phi_j}{\partial \beta} < 0$, and $\frac{\partial \phi_j}{\partial \alpha} > 0$, 
we can obtain by the implicit function theorem 
that 
$\frac{\partial \bar{p}_j}{\partial \beta} = - \frac{\partial \phi_j}{\partial \beta} / \frac{\partial \phi_j}{\partial \bar{p}_j} > 0$, and 
$\frac{\partial \bar{p}_j}{\partial \alpha} = - \frac{\partial \phi_j}{\partial \alpha} / \frac{\partial \phi_j}{\partial \bar{p}_j} < 0$.
\hfill $\Box$

Let $\hat{\xi} := (\underline{\beta}, \overline{\alpha})$. Then we are guaranteed that  
$u_j(p_j, p_{-j}(c_{-j}), c_j, \xi)$ is concave in $p_j$ over
$\mathcal{P}_j(c_j)$.
With the concavity,
conditions (b) and (c) in Assumption \ref{ass:Assumption_1} are satisfied. Without loss of generality, we assume the ambiguity set satisfies condition (d) and that $f_{\theta_j}$ is continuous in $\theta_j$, which satisfies condition (e) as well.
Consequently, we define this interval as firm $j$’s action space. In conclusion, we summarize as follows. 

\begin{theorem}[Existence of BDRNE-MNL]
\label{the:BDRNE-MNL} Let $c \in \mathcal{C}$ be a fixed parameter, and $\bar{p}_j(\cdot, \cdot, \cdot)$ is defined as in Proposition \ref{pro:concavity_u_j}. Then there exists an $n$-tuple
$p^\ast(c):= (p^\ast_1(c_1), \cdots, p^\ast_n(c_n))\in \mathcal{P}(c) := \mathcal{P}_1(c_1) \times \cdots \times \mathcal{P}_n(c_n)$ that satisfies \eqref{eq:BDRNE-MNL}, where $\mathcal{P}_j(c_j):= [c_j, \bar{p}_j(c_j, p^\ast_{-j}(c_{-j}), \hat{\xi})]$, for $j = 1, \cdots, n$.
\end{theorem}


\noindent
\textbf{Proof.}
Under Assumption \ref{ass:Assumption_1}, 
we assert by virtue of 
Theorem \ref{the:existence}, Proposition \ref{pro:concavity_u_j}, and Lemma \ref{lem:barpj} that there exists a 
$p^\ast(c)$ 
satisfying \eqref{eq:BDRNE-MNL}. 
\hfill $\Box$

\subsection{Stability of BDRNE-MNL}

As mentioned earlier, the BDRNE-MNL is dependent on the marginal costs $c$. Therefore, it is natural to discuss the stability of the equilibrium as the marginal cost varies. It is widely known and recognized that when a unique parameterized equilibrium exists within a compact set, the continuity of Nash equilibria is a well-established result, see e.g.~\cite[Proposition 3.2]{bianchi2003note} and \cite{feinstein2022continuity}. The next proposition states the stability of the BDRNE-MNL. 

\begin{proposition}[Stability of BDRNE-MNL]\label{pro:StabilityofBDRNE-MNL} 
Let $p^\ast(c)$ be the unique BDRNE-MNL that satisfies \eqref{eq:BDRNE-MNL}. Then $p^\ast(c)$ is continuous with respect to variations in $c$.
\end{proposition}
\noindent
\textbf{Proof.}
Let $p^\ast(c)$ be the unique BDRNE-MNL that satisfies \eqref{eq:BDRNE-MNL}, for fixed $p^\ast_{-j}(c_{-j})$, we may rewrite the objective of \eqref{eq:BDRNE-MNL} as: 
$$
\hat{u}(p_j, p^\ast_{-j}(c_{-j}), c_j):=
\mathbb{E}_{\theta_j^{N_j}} \left[ \min_{Q\in \mathcal{Q}_{\epsilon_j}^{\theta_j}} \mathbb{E}_{Q} \big[u_j(p_j, p^\ast_{-j}(c_{-j}), c_j, \xi)\big] \right]
- \mathcal{I}_{\mathcal{P}_j(c_j)}(p_j),
$$
where $u_j(p_j, p^\ast_{-j}(c_{-j}), c_j, \xi)$ is concave in $p_j$ over the interval $\mathcal{P}_j(c_j)$ (see Theorem \ref{the:BDRNE-MNL}), and $\mathcal{I}_{\mathcal{P}_j(c_j)}(p_j)$ is the characteristic function defined as
$$
\mathcal{I}_{\mathcal{P}_j(c_j)}(p_j) = 
\begin{cases}
0, & \inmat{for} \;\; p_j\in \mathcal{P}_j(c_j), \\
+\infty, & \inmat{for}\;\;p_j \notin \mathcal{P}_j(c_j).
\end{cases}
$$
The upper level set of $\hat{u}(p_j, p^\ast_{-j}(c_{-j}), c_j)$, denoted as 
$$\mathcal{L}_{\pi}(\hat{u}(p_j, p^\ast_{-j}(c_{-j}), c_j)) := \big\{p_j\mid \hat{u}(p_j, p^\ast_{-j}(c_{-j}), c_j) \geq \pi \big\},$$
is convex and closed for every 
$\pi$.
Thus
$\hat{u}(p_j, p^\ast_{-j}(c_{-j}), c_j)$ is upper semi-continuous and concave.
By \cite[Proposition 3.2]{bianchi2003note}, the BDRNE-MNL is locally Lipschitz in $p_j$, and hence it is stable.
\hfill $\Box$

\section{Numerical results} \label{se:Numerical Results}
In this section, we present specific instances of the BDRNE-MNL model as defined in \eqref{eq:BDRNE-MNL}. We conduct experiments to validate and confirm the accuracy of our theoretical analysis. These test results provide empirical evidence supporting the validity and applicability of the proposed models.

We begin by considering a market involving two mobile phone firms that need to determine the prices of their upcoming new products based solely on their respective marginal costs. In this case, we assume that the performance of the products is measured by a single parameter, which means 
$x$ in \eqref{eq:U_j-MNL} comprises a single component with 
$x_1 = 6$ and $x_2 = 4$, respectively. 
Based on \eqref{eq:U_j-MNL}, we can define the consumer's utility for purchasing each product as follows:
\bgeq
\label{eq:U_j_Numerical}
U_1 = 6\beta - \alpha p_1 + \varepsilon_1,
\quad U_2 = 4\beta - \alpha p_2 + \varepsilon_2,
\edeq
where $\varepsilon_1$ and $\varepsilon_2$ follow Gumbel distributions.
Consequently, we formally formulate the BDRNE-MNL model for the two products as follows:
\begin{equation}\label{eq:Pricing_B1}
\max_{p_1\in \mathcal{P}_1(c_1)} \mathbb{E}_{\theta_1^{N_1}} \left[ \min_{Q\in \mathcal{Q}_{\epsilon_1}^{\theta_1}} \mathbb{E}_{Q} \big[(p_1 - c_1)\frac{e^{6\beta - \alpha p_1}}{1 + e^{6\beta - \alpha p_1} + e^{4\beta - \alpha p_2(c_2)}}\big] \right],
\end{equation}
\begin{equation}\label{eq:Pricing_B2}
\max_{p_2\in \mathcal{P}_2(c_2)} \mathbb{E}_{\theta_2^{N_2}} \left[ \min_{Q\in \mathcal{Q}_{\epsilon_2}^{\theta_2}} \mathbb{E}_{Q} \big[(p_2 - c_2)\frac{e^{4\beta - \alpha p_2}}{1 + e^{4\beta - \alpha p_2} + e^{6\beta - \alpha p_1(c_1)}}\big] \right],
\end{equation}
where $c_1$ and $c_2$ denote the marginal costs of the two products in thousands of dollars. 

In the experiment, we 
follow Shapiro et al.~\cite[Section 4]{shapiro2023bayesian} to consider that both the true distribution and the nominal distribution of the random vector 
$\xi = (\beta, \alpha)$
are joint gamma distributions of
$\Gamma (15, \theta_\beta) \times \Gamma(15, \theta_\alpha)$,
where $\theta_\beta$ and $ \theta_\alpha$ are rate parameters
with $(\theta_\beta^*, \theta_\alpha^*) = (50, 40)$. 
To obtain closed-form posterior updates,
we employ the joint gamma distribution $\Gamma(1, 1) \times \Gamma(1, 1)$, $\Gamma(1, \frac{1}{2}) \times \Gamma(1, \frac{1}{2})$ as the prior distributions of rate parameters 
$(\theta_\beta, \theta_\alpha)$
for firm 1 and firm 2.
Please note that this choice of the conjugate prior distributions is only for computational convenience. 
If the Bayesian updating does not admit a closed-form posterior, we may use Monte Carlo simulation,
such as Markov Chain Monte Carlo (MCMC) methods, to draw samples from the posterior distribution.

We use Algorithm \ref{alg:BDRNE} to solve problems \eqref{eq:Pricing_B1} and \eqref{eq:Pricing_B2} with varying $\epsilon$ and $N$ when $c = (6, 5)$.
Table \ref{table1} shows the performance of firm 1's optimal solution and model value for $\epsilon = 0.01, 0.1, 0.5$, which means $(\epsilon_1, \epsilon_2) = (0.01, 0.01)$, (0.1, 0.1), (0.5, 0.5), and sample size $N = 5$, 20 and 50, respectively. Moreover, we 
compute the equilibrium of problems 
$$\max_{p_j} \mathbb{E}_{\hat{Q}_{N_j}}[u_j(p_j, p_{-j}, c_j, \xi)], \quad \inmat{for}\; j = 1, 2,$$ 
under the empirical distribution $\hat{Q}_{N_j}$.
For further comparative analysis, we consider the DRNE problem, where the ambiguity set is represented by the KL-divergence $\mathcal{Q}_j := 
\big\{ Q: \dd_{KL}\big(Q || \hat{Q}_{N_j}\big) \leq \hat{\epsilon}_j \big\}$ for $\hat{\epsilon} = 0.1$, 
which means $(\hat{\epsilon}_1, \hat{\epsilon}_2) = (0.1, 0.1)$, where $\hat{Q}_{N_j}$ is the nominal distribution based on empirical distribution.
The equilibrium under the true distribution $Q^\ast$ is approximated 
by the empirical equilibrium with $N=200$.

\begin{table}[!ht]\footnotesize
\centering
\begin{tabular}{|c|c|c|c|c|c|c|c|}
\hline
{$N=5$} & $\epsilon = 0.01$ & $\epsilon = 0.1$ & $\epsilon = 0.5$ & $\hat{\epsilon} = 0.1$ & BANE & Empirical & True \\
\hline
optimal solution & 9.2958 & 8.9284 & 8.3611 & 8.9526 & 9.0971 & 10.1887 & 9.5325\\
model value &  0.5862 & 0.4438 & 0.2558 & 0.3582 & 0.6751 & 0.8565 &  0.5416\\
\hline
{$N=20$}  & $\epsilon = 0.01$ & $\epsilon = 0.1$ & $\epsilon = 0.5$ & $\hat{\epsilon} = 0.1$ &  BANE & Empirical & True \\
\hline
optimal solution & 9.4515 &  9.0946 & 8.5538 & 8.6447 & 9.2693 & 9.6250 & 9.5325 \\
model value & 0.4845 & 0.3749 & 0.2279 & 0.2310 & 0.5657 & 0.5788 & 0.5416\\
\hline
{$N=50$} & $\epsilon = 0.01$ & $\epsilon = 0.1$ & $\epsilon = 0.5$ & $\hat{\epsilon} = 0.1$ &  BANE & Empirical & True \\
\hline
optimal solution & 9.6985 & 9.3442 & 8.8048 & 8.6608 & 9.5682 & 9.6373 & 9.5325 \\
model value & 0.5317 & 0.4219 & 0.2711 & 0.2391 & 0.6087 & 0.5949 & 0.5416\\
\hline
\end{tabular}

\captionsetup{font = footnotesize}
\caption{BDRNE-MNL with $\epsilon = 0.01, 0.1, 0.5$, DRNE with $\hat{\epsilon} = 0.1$, BANE, empirical solution and true solution of firm 1 with sample data size $N = 5$, 20, 50, respectively.}
\label{table1}
\end{table}

Figure \ref{fig:P12_N} provides a clearer comparison of solutions from various models. This result is based on a single experiment without averaging over multiple samples. As the sample size $N$ increases, the empirical solution, represented by the blue dash line, gradually stabilizes and converges towards the true solution depicted by the red solid line. As $\epsilon$ decreases, the BDRNE $p_j^\ast$ converges towards the true solution. It is worth noting that BANE is a specific instance of BDRNE with $\epsilon = 0$, its solution is even closer to the true solution than BDRNE with $\epsilon = 0.1, 0.5$. Additionally, BDRNE with $\epsilon = 0.1$ approaches the true solution more closely than DRNE with $\hat{\epsilon} = 0.1$, emphasizing the advantage of BDRNE than DRNE. 

\begin{figure}[!htbp] 
\minipage{0.5\textwidth}
 \centering
  \includegraphics[width=\linewidth]{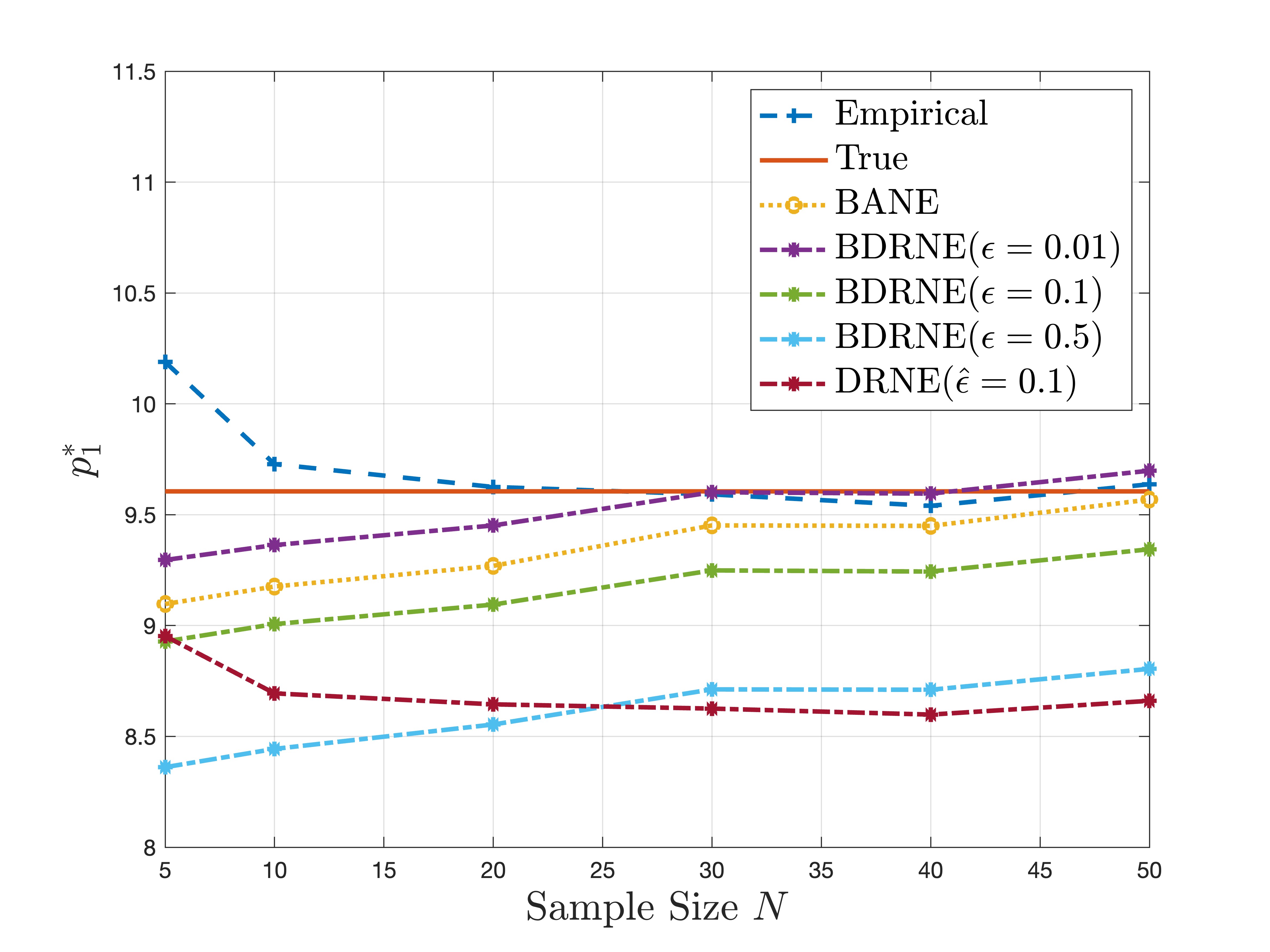}
  \textbf{(a)} 
\endminipage\hfill
\minipage{0.5\textwidth}
  \centering
  \includegraphics[width=\linewidth]{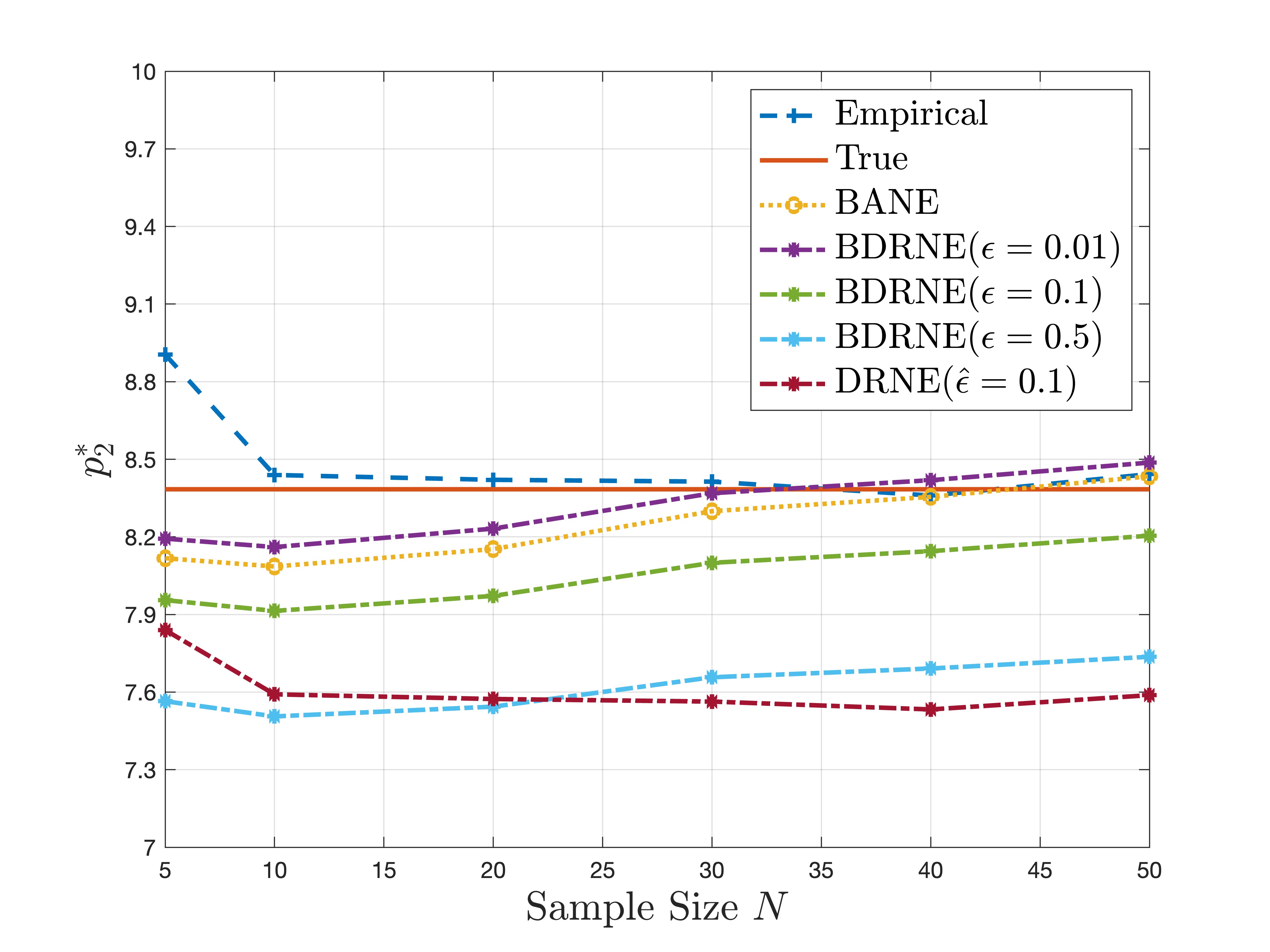}
  \textbf{(b)} 
\endminipage

\captionsetup{font = footnotesize}
\caption{BDRNE-MNL with varying $\epsilon$ values, DRNE with $\hat{\epsilon} = 0.1$, BANE, empirical solution and true solution with increase sample data size, and fixed $c=(6,5)$. \textbf{(a)} $p_1^\ast$ of different models with increase sample size $N$. \textbf{(b)} $p_2^\ast$ of different models with increase sample size $N$.} 
\label{fig:P12_N}
\end{figure}

Figure~\ref{fig:BDRNE_cost} depicts the sensitivity of BDRNE-MNL for problem \eqref{eq:Pricing_B1} and \eqref{eq:Pricing_B2} with varying marginal cost $c_1\in [3, 12]$, while fixed $c_2 = 5$ and sample size $N = 50$. For comparison, we also present the performance of other models such as BANE and DRNE. It is clear that BDRNE exhibits a similar trend to Figure~\ref{fig:P12_N} with varying $\epsilon$ values and performs well than DRNE.
In Figure~\ref{fig:BDRNE_cost} (a), as the marginal cost $c_1$ increases, firm 1 adjusts its price upward to ensure profitability, resulting in a strictly increasing trend for the optimal price function $p_1^\ast$.
In Figure~\ref{fig:BDRNE_cost} (b), with the marginal cost $c_1$ increases, the appeal of product 1 gradually diminishes. Consequently, firm 2 gains more pricing leverage, allowing them to set higher prices. 
The differing slopes of $p_1^\ast$ and $p_2^\ast$ indicate that $p_1^\ast$ is more sensitive to changes in $c_1$.

\begin{figure}[!htbp] 
\minipage{0.5\textwidth}
 \centering
  \includegraphics[width=\linewidth]{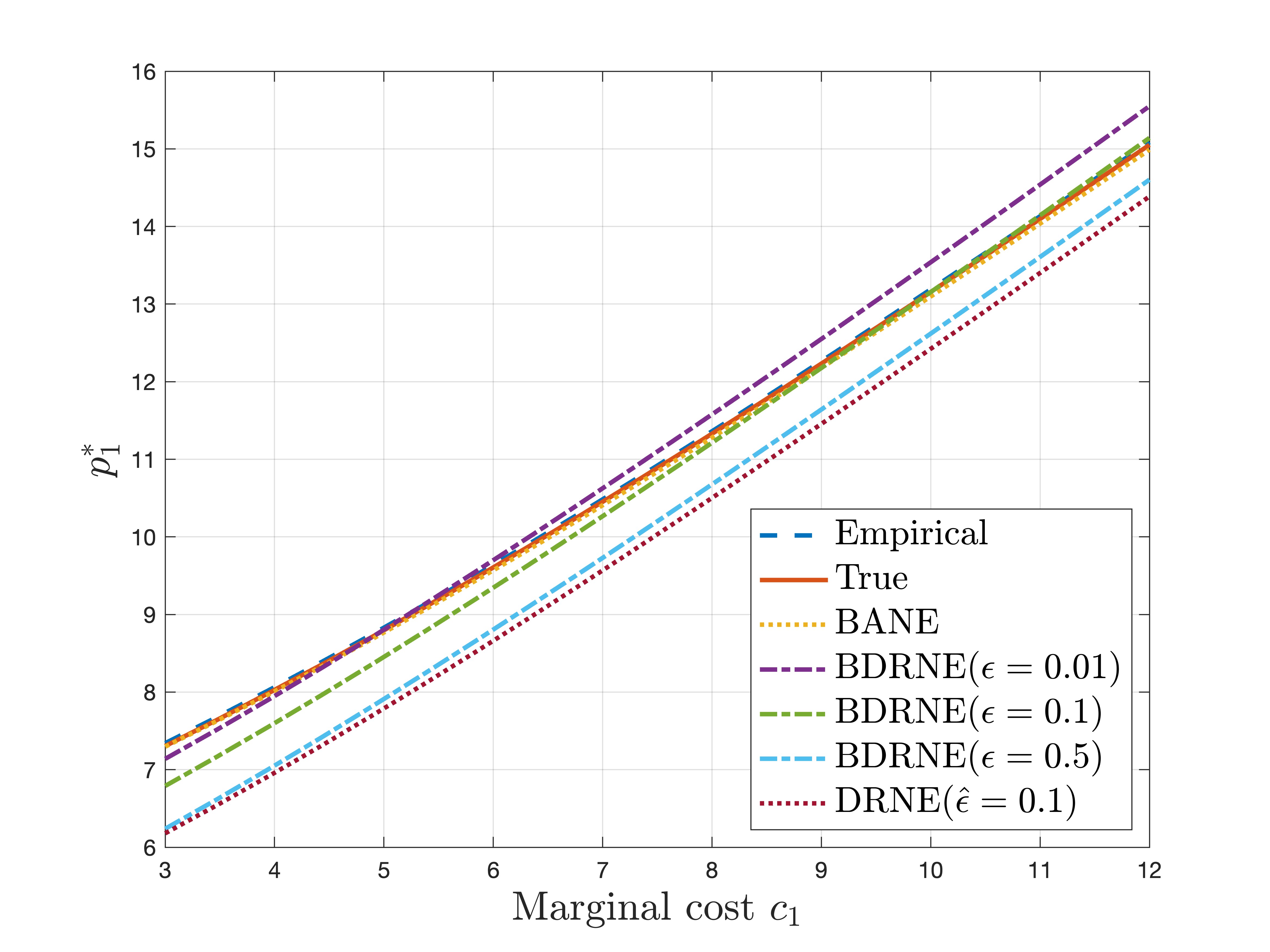}
  \textbf{(a)} 
\endminipage\hfill
\minipage{0.5\textwidth}
  \centering
  \includegraphics[width=\linewidth]{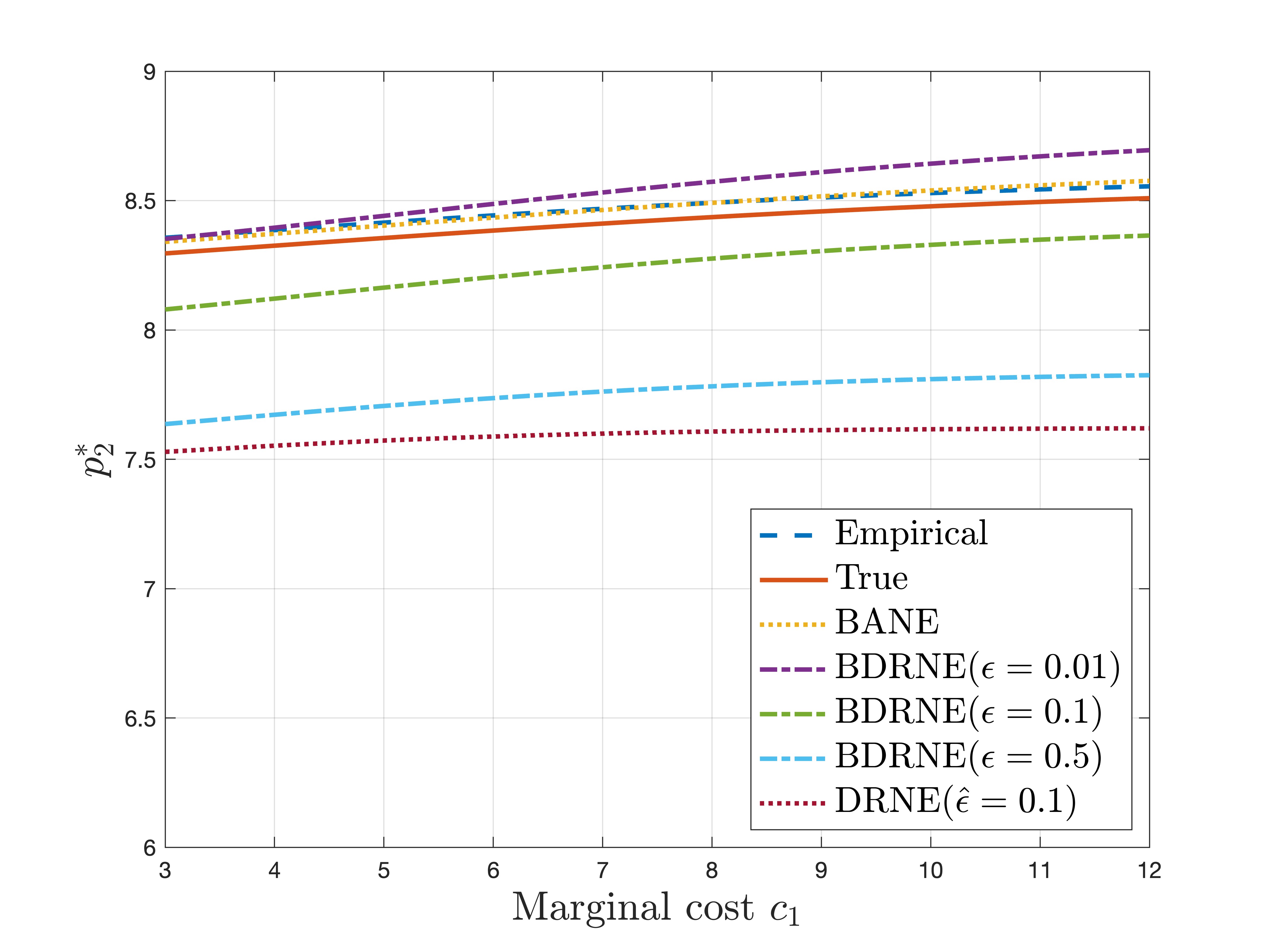}
  \textbf{(b)} 
\endminipage
\captionsetup{font = footnotesize}
\caption{Sensitivity of BDRNE-MNL with varying $\epsilon$ values, DRNE with $\hat{\epsilon} = 0.1$, BANE, empirical solution and true solution with varying marginal cost $c_1$. \textbf{(a)} Sensitivity of $p_1^\ast$ in marginal cost $c_1$.  \textbf{(b)} Sensitivity of $p_2^\ast$ in marginal cost $c_1$.
} 
\label{fig:BDRNE_cost}
\end{figure}

\section{Concluding remarks} \label{se:Concluding remarks}
In this paper, we revisit 
distributionally robust 
Nash equilibrium model by 
applying Bayesian distributionally robust approach recently proposed by Shapiro et al.~\cite{shapiro2023bayesian}. Since the BDRO approach differs significantly from other DRO paradigms in the literature \cite{rahimian2019distributionally}, the new robust equilibrium model may offer an alternative 
modeling approach to analyze the 
outcomes of competition problems.
We have demonstrated how the new equilibrium modeling approach may be applied to 
analyze price competition discrete choice problems.

There are a number of aspects that the new model may be improved or strengthened. First, 
the present way of constructing the ambiguity set may result in large problem size in computing an approximate equilibrium in Algorithm \ref{alg:BDRNE} because the problem size in \eqref{eq:SAA} depends on the posterior distributions of the parameters. This might be addressed by using a new way to define the ambiguity set of the parameters rather than the underlying random variables.
Second, we use Gauss-Seidel-type method to compute an approximate equilibrium in an iterative manner. It is possible to write down the first order optimality conditions of the players' robust optimal decision-making problems and then compute an equilibrium by solving a combined nonlinear complementarity problems when the feasible sets of players have polyhedral structures.
Third, in the asymptotic analysis, we fix the radius of the ambiguity sets. It might be interesting to relate
the radius to the sample size although the subsequent analysis will be rather complicated. 
Fourth, the BDRO model assumes that the samples are i.i.d.~and do not contain noise. In practical data-driven problems, samples may be drawn from empirical data which are not necessarily i.i.d.~\cite{homem2008rates,xu2010uniform}, and they may contain noise \cite{li2023modified}. 
Under these circumstances, it might be interesting to discuss how the equilibrium based on the samples may evolve as the sample size increases and how reliable the equilibrium based on noise data is. 
Fifth, 
as noted by 
Rockafellar \cite{rockafellar1974lagrange,rockafellar1974continuous,rockafellar2019solving},
many decision processes are of a sequential nature and make use of information which is revealed progressively through the observation, at various times,
of random variables with known distributions,
it will be interesting to consider 
multistage BDRNE or BANE.
We leave all these for future research.

\section*{Acknowledgments.}
We would like to thank two referees for their insightful comments which help us significantly strengthen the presentation of this paper. We would also like to thank Editor-in-Chief Professor Liqun Qi for effective handling of the review. Finally, we would like to thank Professor Alex Shapiro for 
valuable discussions during revision of the paper.

\bibliographystyle{abbrv} 
\bibliography{references}

\end{document}